%% file: main.tex
\documentclass[12pt,a4paper,leqno]{article}

\usepackage{amssymb,exscale}
\usepackage[centertags]{amsmath}
\usepackage{amsthm}

\input{preamble}

\begin{document}

\title{From slightly coloured noises\\
  to unitless product systems}
\author{Boris Tsirelson}
\date{}

\maketitle

\stepcounter{footnote}
\footnotetext{Supported in part by the Israel Science Foundation.}

\begin{abstract}
Stationary Gaussian generalized random processes having slowly
decreasing spectral densities give rise to product systems in the
sense of William Arveson (basically, continuous tensor product systems
of Hilbert spaces). A continuum of nonisomorphic unitless product
systems is produced, answering a question of Arveson.
\end{abstract}

\section*{Introduction}
\input{intro}

\section{Measure type spaces and square roots of measures}
\input{sect1}

\section{Gaussian measures, quadratic norms, FHS}
\input{sect2}

\section{Some geometry via measure theory}
\input{sect3}

\section{Density matrices}
\input{sect4}

\section{Fock spaces and tail density matrices}
\input{sect5}

\section{Borel measurability of it all}
\input{sect6}

\section{Sum systems and product systems}
\input{sect7}

\section{The invariant}
\input{sect8}

\section{Slightly coloured noises}
\input{sect9}

\section{The product systems are nonisomorphic and unitless}
\input{sect10}

\bigskip
\filbreak
\begingroup
{
\small
\begin{sc}
\parindent=0pt\baselineskip=12pt

School of Mathematics, Tel Aviv Univ., Tel Aviv
69978, Israel
\emailwww{tsirel@math.tau.ac.il}
{http://www.math.tau.ac.il/$\sim$tsirel/}
\end{sc}
}
\filbreak

\endgroup

\end{document}

%% file: preamble.tex
\numberwithin{equation}{section}

\swapnumbers
\theoremstyle{definition}
\newtheorem{theorem}[equation]{Theorem}
\newtheorem{lemma}[equation]{Lemma}
\newtheorem{proposition}[equation]{Proposition}

\newtheorem{definition}[equation]{Definition}

\renewcommand{\phi}{\varphi}

\renewcommand{\(}{\bigl(}
\renewcommand{\)}{\bigr)\vphantom{)}}

\newcommand{\ip}[2]{\langle#1,#2\rangle}

\newcommand{\Cov}{\operatorname{Cov}}
\newcommand{\Corr}{\operatorname{Corr}}
\newcommand{\Exp}{\operatorname{Exp}}
\newcommand{\Tr}{\operatorname{Tr}}
\newcommand{\mes}{\operatorname{mes}}
\newcommand{\dist}{\operatorname{dist}}
\newcommand{\di}{\operatorname{dim}}
\newcommand{\codim}{\operatorname{codim}}
\newcommand{\spa}{\operatorname{span}}
\newcommand{\closure}{\operatorname{closure}}
\renewcommand{\Re}{\operatorname{Re}}
\renewcommand{\Im}{\operatorname{Im}}
\newcommand{\Ci}{\operatorname{ci}}
\newcommand{\Si}{\operatorname{si}}
\newcommand{\Var}{\operatorname{Var}}
\newcommand{\Const}{\operatorname{Const}}
\newcommand{\const}{\operatorname{const}}
\newcommand{\sgn}{\operatorname{sgn}}

\newcommand{\HS}{\mathcal{HS}}

\newcommand{\bF}{\mathbf F}
\newcommand{\bL}{\mathbf L}
\newcommand{\bA}{\mathbf A}
\newcommand{\Om}{\Omega}
\newcommand{\om}{\omega}
\newcommand{\al}{\alpha}
\newcommand{\be}{\beta}
\newcommand{\eps}{\varepsilon}
\newcommand{\ga}{\gamma}
\newcommand{\Ga}{\Gamma}
\newcommand{\s}{\sigma}
\newcommand{\si}{\sigma}
\newcommand{\de}{\delta}
\newcommand{\E}{\mathcal E}
\newcommand{\F}{\mathcal F}
\newcommand{\A}{\mathcal A}
\newcommand{\B}{\mathcal B}
\newcommand{\D}{\mathcal D}
\newcommand{\U}{\mathcal U}
\newcommand{\cN}{\mathcal N}
\newcommand{\cH}{\mathcal H}
\renewcommand{\P}{\mathcal P}
\newcommand{\la}{\lambda}

\newcommand{\Ex}{\mathbb E}
\newcommand{\R}{\mathbb R}
\newcommand{\Q}{\mathbb Q}
\newcommand{\Z}{\mathbb Z}
\newcommand{\cE}[2]{\mathbb{E}\,\(\,#1\,\big|\,#2\,\)\,}

\newcommand{\sif}{$\sigma$-field}

\newcommand{\II}{\mathcal I}
\newcommand{\I}{\tilde\mathcal I}
\newcommand{\One}{\mathbf1}
\newcommand{\imply}{\;\;\;\Longrightarrow\;\;\;}
\newcommand{\imp}{$ \Longrightarrow $ }
\renewcommand{\equiv}{\;\;\;\Longleftrightarrow\;\;\;}

\def\emailwww#1#2{\par\qquad {\tt #1}\par\qquad {\tt #2}\medskip}

%% file: intro.tex
The white noise is a Gaussian stationary generalized random process
whose restrictions to adjacent intervals $ (a,b) $ and $ (b,c) $ are
independent. In contrast, for a continuous process, its restrictions
to $ (a,b) $ and $ (b,c) $ are heavily dependent via the value at $ b
$. Such a dependence cannot be described by a probability density; the
joint distribution is singular w.r.t.\ the product of marginal
distributions. By a slightly coloured noise I mean a Gaussian
stationary generalized random process such that the distribution of
its restriction to $ (a,c) $ is absolutely continuous w.r.t.\ the
product of the distributions of its restrictions to $ (a,b) $ and $
(b,c) $, whenever $ -\infty < a < b < c < +\infty $.

In the Hilbert space of all linear functionals of the white noise,
every interval $ (a,b) $ determines a subspace $ G_{a,b} $ satisfying $
G_{a,b} \oplus G_{b,c} = G_{a,c} $; the whole space is a direct
integral (a continuous direct sum). For arbitrary (not just linear)
functionals the relation is multiplicative rather than additive:
\[
H_{a,b} \otimes H_{b,c} = H_{a,c} \, ;
\]
the whole space is a continuous tensor product, in other words, a
product system. The constant $ 1 $ may be treated as an element $
\One_{a,b} \in H_{a,b} $, satisfying
\[
\One_{a,b} \otimes \One_{b,c} = \One_{a,c} \, ;
\]
such a multiplicative family is called a unit (of a product
system). There exist product systems with many units, with a single
unit (up to a natural equivalence), and unitless (with no
unit). However, the theory of unitless product systems suffers from
lack of rich sources of examples. Slightly coloured noises are such a
source, rich enough for producing a continuum of nonisomorphic
unitless product systems.

\begin{quote}
``We believe that there should be a natural way of constructing such
product systems, and we offer that as a basic unsolved problem. The
fact that we do not yet know how to solve it shows how poorly
understood continuous tensor products are today.''\\
\hspace*{\fill} Arveson 1994 \cite[p.~5]{Ar94}.
\end{quote}

After producing a continuum of nonisomorphic product systems with
units \cite{TsAr} I was asked by Arveson (private communication,
January 2000) about a continuum of nonisomorphic unitless product
systems. His question is answered here by using a construction that
was outlined by Tsirelson and Vershik \cite[Sect.~1c]{TV} with no
proofs.

\begin{quote}
``This example may be considered as a commutative (bosonic)
counterpart of Power's noncommutative (fermionic) example of a
non-Fock factorization over $ \R $.''\\
\hspace*{\fill} Tsirelson and Vershik 1998 \cite[p.~91]{TV}.
\end{quote}

\begin{center}
\textsc{Probabilistic aspects}
\end{center}

Probabilistic results are summarized here in probabilistic language
(while the rest of the paper is written in rather analytical
language).

Consider a Gaussian stationary generalized random process $
(\xi_t)_{t\in\R} $ whose covariation function $ B(t) = \Cov ( \xi_s,
\xi_{s+t} ) $ is positive, decreasing and convex on $ (0,\infty) $,
and
\[
B(t) = \frac1{ |t| \ln^\al (1/|t|) } \quad \text{for all $ t $ small
enough} \, ;
\]
here $ \al \in (1,\infty) $ is a parameter.

Then the joint distribution of random variables
\[
X_k = \int_0^{2\pi} e^{ikt} \xi_t \, dt \qquad ( k =
\dotsc,-2,-1,0,1,2,\dotsc )
\]
has a density (finite and strictly positive almost everywhere) w.r.t.\
the product of corresponding one-dimensional distributions $ N ( \Ex
X_k, \Var X_k ) $.

The same holds for a larger family $
(\dotsc,X_{-1},Y_{-1},X_0,Y_0,X_1,Y_1,X_2,Y_2,\dotsc) $
of random variables, where $ X_k $ are as before, and $ Y_k =
\int_{-2\pi}^0 e^{ikt} \xi_t \, dt $.

Another result. Take some $ \eps_n \in (0,1) $ and consider random
variables
\[
Z_n = \frac1{\eps_n} \sum_{k=0}^{n-1} \int_{k/n}^{(k+\eps_n)/n} \xi_t
\, dt \, , \qquad Z = \int_0^1 \xi_t \, dt \, .
\]

If $ \eps_n \ln^{\al-1} n \to \infty $ then $ Z_n \to Z $ in $ L_2 $.

If $ \eps_n \ln^{\al-1} n \to 0 $ then $ \Var(Z_n) \to \infty $ and
the correlation coefficient $ \Corr (Z_n,Z) \to 0 $.

\begin{sloppypar}
If $ \lim_n ( \eps_n \ln^{\al-1} n ) \in (0,\infty) $ then $ \lim_n
\Var(Z_n) \in (0,\infty) $ and $ \lim_n \Corr (Z_n,Z) \in (0,1) $.
\end{sloppypar}

For detail see Sect.~10. The reader interested just in these
probabilistic statements may skip sections 1--9 in which case,
however, he/she should bypass some points of analytical nature in
Sect.~10, and restore some proofs omitted since they are not necessary
from the analytical viewpoint.

\begin{center}
\textsc{Quantal aspects}
\end{center}

Acquaintance with quantum theory is not needed for reading the rest of
the paper, but should help to understand the idea as explained here.

A product system may be thought of as a local quantum field over the
one-dimensional space $ \R $ (just space, no time at all; see also
Arveson \cite{Ar99} for a better, dynamical interpretation of product
systems). The whole field is a quantum system, and its restriction to
$ (0,1) $ is a subsystem. A unit vector of $ H_{0,1} $ describes a
pure state of the subsystem (though in general the subsystem and the
rest of the system are entangled).

Local Hilbert spaces (and algebras) are ascribed to intervals, as well
as to more general regions, consisting of a finite number of
intervals. Introduce a region $ E_n $ consisting of $ n $ small
equidistant intervals of equal length,
\[
E_n = \Big( 0, \frac{\eps_n}{n} \Big) \cup \Big( \frac1n,
\frac{1+\eps_n}{n} \Big) \cup \dots \cup \Big( \frac{n-1}n,
\frac{n-1+\eps_n}{n} \Big)
\]
and consider the corresponding subsystem, described by its Hilbert
space $ H_{E_n} $. For a fixed $ n $, the subsystem may be entangled
or not (with the rest). If the product system has a unit then there is
a `white state' with no spatial correlations. It makes $ H_{E_n} $
disentangled for all $ n $ simultaneously. That is the case for a
product system constructed out of the white noise.

A unitless product system is constructed out of the slightly coloured
noise mentioned in `Probabilistic aspects'. Spatial correlations
inherent to the noise are inherited by quantum states. Subsystems $
H_{E_n} $ can be disentangled for finitely many $ n $, but not for all
$ n $ simultaneously.

\begin{sloppypar}
Asymptotic behavior of subsystems $ H_{E_n} $ depends crucially on $
\lim ( \eps_n \ln^{\al-1} n ) $, as is suggested by properties of
random variables $ Z_n $ (see `Probabilistic aspects'). If $ \lim (
\eps_n \ln^{\al-1} n ) = 0 $ then $ E_n $ are a tail sequence in the
sense that the mixed state of the subsystem (the density matrix on $
H_{E_n} $) has a universal asymptotics, irrespective of the state of
the system. The situation is different if $ \lim ( \eps_n \ln^{\al-1}
n ) \ne 0 $. This is why product systems for different $ \al $ are
nonisomorphic.
\end{sloppypar}

\begin{center}
\textsc{Aspects of functional analysis}
\end{center}

A classical construction going back to Fock may be outlined as
follows. One starts with a Hilbert space $ G $. One identifies $ G
$ with the space of all measurable linear functionals over a Gaussian
measure $ \ga $. One gets another Hilbert space $ H = L_2 (\ga) $ that
may be denoted $ H = \Exp G $, since $ G = G_1 \oplus G_2 $ implies $
H = H_1 \otimes H_2 $ (via $ \ga = \ga_1 \otimes \ga_2 $). A
continuous direct sum, $ G = \int^\oplus G_\zeta \, d\zeta $
(basically the same as a projection-valued measure well-known in
spectral theory) leads to a continuous tensor product. Classical
product systems, obtained this way, contain units.

My modification of the classical construction is rather innocent (but
surprisingly powerful). Basically, the subspaces $ G_1, G_2 $ are
allowed to be slightly nonorthogonal, without destroying the relation
$ H = H_1 \otimes H_2 $. In terms of Gaussian measures, the relation $
\ga = \ga_1 \otimes \ga_2 $ is generalized to $ \ga = p \cdot ( \ga_1
\otimes \ga_2 ) $ where $ p $ is a density (that is, Radon-Nikodym
derivative). The density $ p $ is inserted properly into the formula
for $ h_1 \otimes h_2 \in H $.

The classical language (of Hilbert spaces, Gaussian measures, spaces $
L_2(\ga) $ etc.) is not well-suited to the modified construction. A
modified language, used in the paper, stipulates larger invariance
(symmetry) groups. Namely, the `Hilbert space' structure corresponds
to the group of all unitary operators, $ U^* U = I $. I use a larger
group consisting of all invertible $ U $ such that $ U^* U - I $ is a
Hilbert-Schmidt operator. (These $ U $ are called `equivalence
operators' in Feldman's well-known paper \cite{Fe} on equivalence of
Gaussian measures). The corresponding structure, weaker than `Hilbert
space' structure but stronger than `linear topological space'
structure, is defined in Sect.~2 under the name `FHS-space'. An
equivalence class (in Feldman's sense) of norms is used rather than a
single norm.

Accordingly, a measure space $ (\Om,\F,P) $ is replaced with a
`measure type space' $ (\Om,\F,\P) $. An equivalence class $ \P $ of
measures is used rather than a single measure $ P $. See Sect.~1. It
appears that $ L_2 (\Om,\F,\P) $ can be defined naturally (in addition
to the usual $ L_2 (\Om,\F,P) $); see Sect.~1.

Also, a `Gaussian type space' defined in Sect.~2 stipulates an
equivalence class of Gaussian measures rather than a single Gaussian
measure.

For a technical reason we need Borel measurability of several natural
constructions. For example, the orthogonal projection of a vector to a
subspace (in a Hilbert space) is a jointly Borel measurable function
of the point and the subspace, provided that the set of all subspaces
is equipped with its natural Borel structure. Similarly, the
conditional expectation is a jointly Borel measurable function of a
random variable and a sub-\sif. See Sect.~6 for detail.

%% file: sect1.tex
Let $ \Om $ be a nonempty set, $ \F $ a \sif\ of its subsets, and $
\mu, \nu $ measures\footnote{%
 I assume always, that every measure space $ (\Om,\F,\mu) $
 is a Lebesgue-Rokhlin space (that is, isomorphic $ \bmod 0 $ to
 an interval with Lebesgue measure, to a finite or countable set of
 atoms, or a combination of both). However, the assumption is not
 really used in this work. We deal with such objects as $ L_2
 (\Om,\F,\mu) $; the latter must be separable; other properties of $
 (\Om,\F,\mu) $ do not matter. I assume also that all considered
 \sif s ($ \F $ itself, and its sub-\sif s) contain all negligible
 sets. Everything is treated $ \bmod 0 $, that is, up to negligible
 sets.}
(real-valued, finite, positive) on $ \F $. One
says that $ \mu, \nu $ are \emph{equivalent,} if they are mutually
absolutely continuous. Let $ \P $ be an equivalence class of
probability measures on $ (\Om,\F) $. That is, every $ P \in \P $ is a
measure on $ (\Om,\F) $ satisfying $ P(\Om) = 1 $, and every $ P_1,
P_2 \in \P $ are equivalent, and $ \P $ contains every probability
measure equivalent to a measure of $ \P $ (and of course, $ \P $ is
nonempty). One says that $ \P $ is a type of measure, and $
(\Om,\F,\P) $ is a \emph{measure type space.}

\begin{sloppypar}
The linear topological (metrizable, but not locally convex) space $
L_0 (\Om,\F,\P) $ consists of all (equivalence classes of) measurable
functions on $ \Om $. Its topology corresponds to convergence in
measure (in probability), irrespective of the choice of a measure $ P
\in \P $.
\end{sloppypar}

Hilbert spaces $ L_2 (\Om,\F,P_1) $ and $ L_2 (\Om,\F,P_2) $ for $
P_1, P_2 \in \P $ differ (unless $ \exists \eps \;\; \eps P_1 \le P_2
\le (1/\eps) P_1 $). However, they are in a \emph{natural} unitary
correspondence. Namely, $ \psi_1 \in L_2 (\Om,\F,P_1) $ corresponds to
$ \psi_2 \in L_2 (\Om,\F,P_2) $ when $ \psi_2 = \sqrt{ \frac{P_1}{P_2}
} \psi_1 $; here and henceforth $ \frac{P_1}{P_2} $ stands for the
Radon-Nikodym density, and I write $ P_1 = \frac{P_1}{P_2} \cdot P_2
$. We'll glue all $ L_2 (\Om,\F,P) $ together, forming $ L_2
(\Om,\F,\P) $ that contains $ \sqrt P $ for all $ P \in \P $, as we'll
see soon.

There are several reasons for introducing `square roots of
measures'. One reason. In quantum mechanics, Schr\"odinger's wave
function $ \psi(x) $ determines a probability distribution $
|\psi(x)|^2 \, dx $. However, for infinitely many degrees of freedom
we have no Lebesgue measure ($ dx $). It could be convenient to
describe a quantum state by an object that combines a measure and
phases, something like $ \psi(x) \sqrt{dx} $.

Old-fashioned tensor analysis stipulates a notion of a relative
tensor, in particular, a relative scalar of a given weight (see for
instance \cite[item 156 on pp.~345--346]{Br}). A density of a measure
on a manifold is a relative scalar field of weight $ 1 $. Relative
scalar fields of weight $ 1/2 $ are smooth finite-dimensional `square
roots of measures' considered below.

Define $ L_2 (\Om,\F,\P) $, denoted also by $ L_2 (\P) $ for short, as
the Hilbert space of all families $ \psi = ( \psi_P )_{P\in\P} $ such
that $ \psi_P \in L_2 (\Om,\F,P) $ for every $ P \in \P $, and
\begin{equation}\label{11a}
\psi_{P_2} = \sqrt{ \frac{P_1}{P_2} } \psi_{P_1}
\end{equation}
for all $ P_1,P_2 \in \P $. Linear operations and scalar product are
defined naturally; for every $ P \in \P $, the map $ L_2 (\P) \ni
\psi \mapsto \psi_P \in L_2 (P) $ is unitary.

\begin{sloppypar}
Given $ P \in \P $, we define $ \sqrt P \in L_2 (\P) $ by $ (\sqrt
P)_P = 1 $ (identically on $ \Om $). More generally, given $ P \in \P
$ and $ f \in
L_2 (P) $, we define $ f \sqrt P \in L_2 (\P) $ by $ (f\sqrt P)_P = f
$. Thus, $ \psi = \psi_P \sqrt P $ for all $ P \in \P $, $ \psi \in
L_2 (\P) $. Now we may replace the notation $ \psi_P $ by a more
expressive notation
\[
\psi_P = \frac{ \psi }{ \sqrt P } \, .
\]
Given $ \psi', \psi'' \in L_2 (\P) $, their scalar product is $
\langle \psi', \psi'' \rangle = \int \frac{\psi'}{\sqrt P}
\frac{\psi''}{\sqrt P} \, dP $; the latter does not depend on $ P \in
\P $. Heuristically we could write $ \langle \psi', \psi'' \rangle =
\int \psi'(\om) \psi''(\om) $, however, I prefer the notation
\[
\langle \psi', \psi'' \rangle = \int \frac{\psi'(\om)}{\sqrt{P(d\om)}}
\frac{\psi''(\om)}{\sqrt{P(d\om)}} \, P(d\om) = \int \psi'(\om)
\psi''(\om) \frac{d\om}{d\om} \, ;
\]
being rather awkward and illogical,\footnote{%
 Neither $ \psi(\om) $ nor $ \psi(d\om) $ is a logical notation for an
 object of the form $ f(\om) \sqrt{P(d\om)} $.}
it still helps.
Given $ \psi', \psi'' \in L_2 (\P) $, we define the signed measure $
\psi' \psi'' $ by $ \frac{\psi'\psi''}{P} = \frac{\psi'}{\sqrt P}
\frac{\psi''}{\sqrt P} $, thus $ \ip{\psi'}{\psi''} = (
\psi' \psi'' ) (\Om) $, the measure of the whole space.
Note also that $ \langle \sqrt{P_1}, \sqrt{P_2} \rangle = \int \sqrt{
P_1(d\om) } \sqrt{ P_2(d\om) } $ was used in Kakutani's well-known
work \cite{Ka} about equivalence of product measures. Two natural
metrics on $ \P $ define the same topology (the only one used here); I
mean the variation distance, $ \| P_1 - P_2 \| = \int \big|
\frac{P_1}P - \frac{P_2}P \big| \, dP $, and the angle in $ L_2 (\P)
$, $ \arccos \ip{\sqrt{P_1}}{\sqrt{P_2}} $; the proof is left to the
reader.\footnote{%
 In fact, $ \( \frac12 \| P_1 - P_2 \| \)^2 \le 1 - \ip
 {\sqrt{P_1}}{\sqrt{P_2}} \le \frac12 \| P_1 - P_2 \| $, which
 is not used here.}
\end{sloppypar}

%% file: sect2.tex
One may introduce a Gaussian measure as a probability measure on a
Banach (or Hilbert, or Frechet, etc) space, such that every
(continuous) linear functional has a normal distribution. Such a
viewpoint is convenient for heuristic thinking but not for a formal
presentation, since topological structures on the linear space with
measure are in fact irrelevant. Following an old advice of A.~Vershik,
I prefer discarding any topological structure on the (Banach, etc.)
space and even keeping implicit its linear structure. Everything can
be formulated in terms of two different spaces, a probability space
and a Hilbert space; the latter plays the role of the space of all
measurable linear functionals, as well as its dual, the space of all
admissible shifts.

Let $ (\Om,\F,P) $ be a probability space and $ G \subset L_2
(\Om,\F,P) $ a (closed) linear subspace, containing constants and
generating the whole \sif\ $ \F $. We call $ G $ a \emph{Gaussian
space,} if every $ g \in G $ has a normal distribution; the latter is
therefore $ N ( \Ex g, \Var g ) = N \( \int g \, dP , \int g^2 \, dP -
(\int g \, dP)^2 \) $. Up to isomorphism, there is exactly one
Gaussian space in every dimension ($ 0,1,2,\dots $ or $ \infty
$). That is, if $ G_1 \subset L_2 (\Om_1,\F_1,P_1) $ and $ G_2 \subset
L_2 (\Om_2,\F_2,P_2) $ are Gaussian spaces and $ \dim G_1 = \dim G_2
$, then there exists an isomorphism (invertible measure preserving
map) between the probability spaces $ (\Om_1,\F_1,P_1) $ and $
(\Om_2,\F_2,P_2) $ that induces an isometry betwen $ G_1 $ and $ G_2
$. Moreover, each isometry between $ G_1 $ and $ G_2 $ corresponds to
exactly one isomorphism of probability spaces.

The standard model is the space $ \R^\infty = \R \times \R \times
\dots $ of all sequences of real numbers, equipped with the product
measure $ \ga^\infty = \ga^1 \otimes \ga^1 \otimes \dots $ where $
\ga^1 $ is the standard normal distribution $ N(0,1) $. Here $ G $
consists of all measurable linear functionals $ \R^\infty \ni (x_1,
x_2, \dots) \mapsto \sum c_k x_k \in \R $ with $ \sum c_k^2 < \infty
$; that is, $ l_2 $ is the standard model of $ G $. Any other model is
necessarily isomorphic to the standard model, as far as $ G $ is of
infinite dimension; straightforward modifications for a finite
dimension are left to the reader.

The same $ G $ over $ \R^\infty $ is also a Gaussian space w.r.t.\
many other measures on $ \R^\infty $ equivalent to the standard
measure $ \ga^\infty $. In particular, consider the product measure $
\ga $ on $ \R^\infty $ whose $ k $-th factor is the normal
distribution $ N (m_k, \s_k^2) $ with given parameters $ m_k \in
(-\infty, +\infty) $, $ \s_k \in (0,+\infty) $. The well-known S.\
Kakutani's theorem on equivalence of infinite product measures
\cite{Ka} shows that $ \ga $ is equivalent to $ \ga^\infty $ if and
only if
\[
\sum_k (\s_k-1)^2 < \infty \quad \text{and} \quad \sum_k m_k^2 <
\infty \, .
\]
Replacing coordinate axes by arbitrary orthogonal (in $ l_2 $) axes
one gets all Gaussian measures on $ \R^\infty $ equivalent to $
\ga^\infty $, which is basically the well-known criterion (Feldman
\cite{Fe}, Hajek \cite{Ha}, Segal \cite{Se}).

In terms of a Gaussian space $ G \subset L_2 (\Om,\F,P) $ we introduce
the set $ \P $ of all probability measures on $ (\Om,\F) $ equivalent
to $ P $, and its subset $ \P_G $ consisting of all $ \ga \in \P $
such that $ G $ is also a Gaussian space in $ L_2 (\Om,\F,\ga) $;
these measures will be called \emph{Gaussian measures w.r.t.\ $ G $.}
Now we forget the initial measure $ P $; each measure of $ \P_G $ may
serve as $ P $ equally well. Up to isomorphism, the structure $
(\Om,\F,\P,G) $ is uniquely determined by $ \dim G $. I call $
(\Om,\F,\P,G) $ a \emph{Gaussian type space.}

We introduce a quotient space $ G_0 = G / \Const $ (where $ \Const $
is the one-dimensional space of constants) and its dual space $ G^0 $
(consisting of all continuous linear functionals on $ G_0 $). Every $
\ga \in \Ga $ determines a Hilbert norm $ \|\cdot\|_\ga $ on $ G_0 $
(that is, $ (G_0, \|\cdot\|_\ga) $ is a Hilbert space) via
\[
\| g \|^2_\ga = \Var_\ga (g) = \int g^2 \, d\ga - \Big( \int g \, d\ga
\Big)^2 \, ,
\]
which is insensitive to any constant added to $ g \in G $. Each norm $
\|\cdot\|_\ga $ defines the same topology on $ G_0 $. However, the
norms are equivalent in a stronger sense, namely, for any $ \ga_1,
\ga_2 \in \Ga $ there exists a basis $ (g_k) $ of $ G_0 $, orthogonal
for both norms and such that the numbers $ \la_k = \| g_k \|_{\ga_2} /
\| g_k \|_{\ga_1} $ satisfy $ \la_k > 0 $ and $ \sum (\la_k-1)^2 <
\infty $ (which can be reformulated in terms of Hilbert-Schmidt
operators; namely, the unit operator $ (G_0,\|\cdot\|_{\ga_1}) \ni g
 \mapsto g \in (G_0,\|\cdot\|_{\ga_2}) $ must be an equivalence
 operator, as defined by Feldman \cite[Def.~1]{Fe}). I'll call such
 norms \emph{FHS-equivalent.}\footnote{%
 You may interpret FHS as Feldman-Hajek-Segal, or alternatively as
 Feldman-Hilbert-Schmidt.}
Of course, it is an equivalence relation \cite[Lemma 2]{Fe}. Numbers $
 \la_k $ determine the distance between $ \ga_1, \ga_2 $
(provided that $ \ga_1, \ga_2 $ have the same mean, that is, $ \int g
 \, d\ga_1 = \int g \, d\ga_2 $ for all $ g \in G $); namely, a simple
(basically, one-dimensional) calculation gives
\[
\ip{ \sqrt{\ga_1} }{ \sqrt{\ga_2} } = \prod_k \bigg( \frac{
\la_k^{-1/2} + \la_k^{1/2} }{ 2 } \bigg)^{-1/2} \, ;
\]
of course, convergence of the product is equivalent to convergence of
the series $ \sum_k (\la_k-1)^2 $.

\begin{definition}\label{2.1}
An \emph{FHS-space} is a pair $ (H,\cN) $ of a linear space $ H $ and
a set $ \cN $ of norms on $ H $ such that

(a) every norm of $ \cN $ turns $ H $ into a separable Hilbert space;

(b) all norms of $ \cN $ are pairwise FHS-equivalent;

(c) every norm FHS-equivalent to a norm of $ \cN $ belongs to $ \cN
$.

\noindent Norms belonging to $ \cN $ will be called \emph{admissible
norms} on $ H $.
\end{definition}

Isomorphisms of FHS-spaces are basically the same as Feldman's
equivalence operators \cite[Def.~1]{Fe}.
Up to isomorphism, there is exactly one FHS-space in every dimension
($ 0,1,2,\dots $ or $ \infty $); just the same situation as for
(separable) Hilbert spaces.

Every separable Hilbert space is isometric to a Gaussian space (over
some probability space), and this superstructure brings no
arbitrariness, as far as it is considered up to isomorphism. (One
speaks about \emph{the} isonormal random process on any given Hilbert
space.) Similarly, every FHS-space may be identified with $ G / \Const
$ for \emph{the} Gaussian type space $ (\Om,\F,\P,G) $ of the
corresponding dimension. Thus, in principle one may prove a purely
geometric statements about an FHS-space via Gaussian measures, and
this way is indeed used in the next section. The
correspondence between $ \cN $ (admissible norms) and $ \P_G $
(Gaussian measures) transfers to $ \cN $ the natural topology of $
\P_G $ (inherited from $ \P $; recall the end of Sect.~1). In terms of
the numbers $ \la_k $, a basis of neighborhoods of a given admissible
norm may be written as $ \sum_k (\la_k-1)^2 < \eps $.

The (additive group of) space $ G^0 $, dual to $ G_0 $, acts on $
(\Om,\F,\P) $ by automorphisms (invertible measurable transformations
preserving $ \P $) $ U_x $, $ x \in G^0 $, such that
\[
g (U_x \om) - g (\om) = \langle g, x \rangle \quad \text{for almost
all } \om \in \Om
\]
for all $ g \in G $. Clearly, $ U_x $ is uniquely
determined. Existence of these $ U_x $ may be checked just for the
standard model, which makes it evident: $ G^0 = l_2 $ acts on $
\R^\infty $ by shifts, $ U_x (\om) = \om + x $. The map $ U_x $ sends
each Gaussian measure $ \ga_1 \in \P_G $ into another Gaussian measure
$ \ga_2 \in \P_G $ such that $ \|\cdot\|_{\ga_1} = \|\cdot\|_{\ga_2} $
on $ G_0 $, and $ \int g \, d\ga_2 - \int g \, d\ga_1 = \langle g, x
\rangle $ for all $ g \in G $. A simple (basically, one-dimensional)
calculation gives
\[
\ip{ \sqrt{\ga_1} }{ \sqrt{\ga_2} } = \exp \( -\tfrac18 \| x \|^2 \) \,
,
\]
where $ \| \cdot \| = \| \cdot \|_{\ga_1}  = \| \cdot \|_{\ga_2} $. In
fact, $ \sqrt{ \ga_1 } \cdot \sqrt{ \ga_2 } = \exp \( -\tfrac18 \| x
\|^2 \) \ga $, where $ \ga $ is the image of $ \ga_1 $ under $ U_{x/2}
$.

Turn to the Hilbert space $ L_2 (\Om,\F,\P) $ of `square roots of
measures'. Transformations $ U_x $ induce unitary operators on $ L_2
(\P) $; I denote them by $ U_x $, too. Namely, for every $ \psi \in
L_2(\P) $ and $ P \in \P $,
\[
\frac{ \psi }{ \sqrt P } (\om) = \frac{ U_x \psi }{ \sqrt{U_x P} }
(U_x \om) \, ,
\]
where $ U_x P $ (denoted also by $ P \circ U_x^{-1} $) is the image of
$ P $ under $ U_x $. Thus, the FHS-space $ G^0 $ acts unitarily on the
Hilbert space $ L_2 (\P) $.

There is also a natural projective action of the other FHS-space, $
G_0 $, on $ L_2 (\P) $. Namely, every $ g \in G $ determines a unitary
operator $ V_g : L_2 (\P) \to L_2 (\P) $,
\[
V_g \psi = e^{ig} \psi \, ,
\]
the multiplication by the function $ \om \mapsto e^{ig(\om)} $. Given
$ y \in G_0 $, we get $ V_y $ determined up to a phase factor, which
means a \emph{projective} action. In fact, one can choose phase
factors getting a unitary representation (which is evident for the
standard model). Anyway,
\[
V_y U_x = e^{i \langle x, y \rangle} U_x V_y \, ,
\]
the well-known Weyl form of Canonical Commutation Relations. In this
context, vectors $ \psi \in L_2 (\P) $ of the form $ \psi = V_y
\sqrt\ga $, $ \ga \in \Ga $, $ y \in G_0 $, are known as coherent
states, or quasi-free pure states, or Gaussian pure states. It is easy
to see that
\[
\ip{ V_y \sqrt\ga }{ \sqrt\ga } = \exp \( -\tfrac12 \| y \|_\ga^2 \)
\]
for $ y \in G_0 $.

%% file: sect3.tex
Many statements of this section are purely geometric and probably
could be proved within the Hilbert space geometry, but I find it
easier to prove them via measure theory. Usually, the corresponding
properties of Gaussian measures hold for arbitrary (non-Gaussian)
measures as well.

The formula $ H = H_1 \oplus H_2 $ has several interpretations. It
always implies that $ H_1, H_2 $ are subspaces of $ H $ and every
vector $ h \in H $ has a unique representation as $ h_1 + h_2 $ where
$ h_1 \in H_1 $, $ h_2 \in H_2 $. However, when $ H $ is a Hilbert
space, one usually stipulates that $ H_1, H_2 $ are orthogonal. When $
H $ is an FHS-space, we treat
\[
H = H_1 \oplus H_2
\]
as follows: there exists an \emph{admissible}\footnote{%
 Recall Definition \ref{2.1}.}
norm on $ H $ that makes $ H_1, H_2 $ orthogonal (and of course, they
span the whole $ H $). Similarly, $ H = H_1 \oplus \dots \oplus H_n $
means existence of an admissible norm that makes $ H_1, \dots, H_n $
orthogonal (also, they span $ H $).

\begin{proposition}\label{3.1}
Let $ H $ be an FHS-space, and $ H_1, H_{12}, H_{23}, H_3 $ its
subspaces such that $ H_1 \subset H_{12} $, $ H_{23} \supset H_3 $,
and
\[
H_{12} \oplus H_3 = H = H_1 \oplus H_{23} \, .
\]
Then the subspace $ H_2 = H_{12} \cap H_{23} $ satisfies
\[
H_1 \oplus H_2 \oplus H_3 = H \, ,
H_1 \oplus H_2 = H_{12} \, ,
H_2 \oplus H_3 = H_{23} \, .
\]
\end{proposition}

Note. Given that $ H_1 \oplus H_2 \oplus H_3 = H $, other relations $
H_1 \oplus H_2 = H_{12} $, $ H_2 \oplus H_3 = H_{23} $ may be treated
in the topological sense.

\begin{proof}
The decomposition $ H = H_1 \oplus H_{23} $ determines a projection $
P_1 : H \to H $ such that $ P_1 H = H_1 $ and $ 1-P_1 = P_{23} $ is
also a projection, $ P_{23} H =H_{23} $. The same for $ P_3 $ and $
P_{12} = 1 - P_3 $. The inclusion $ H_1 \subset H_{12} $ gives $ P_1 =
P_{12} P_1 = (1-P_3) P_1 $, that is, $ P_3 P_1 = 0 $; similarly, $ P_1
P_3 = 0 $. So, our projections commute with each other. Introducing
\[
P_2 = P_{12} P_{23} = P_{23} P_{12} = (1-P_1) (1-P_3) = 1 - P_1 - P_3
\, ,
\]
we have $ P_2^2 = P_{12} P_{23} P_{12} P_{23} = P_{12}^2 P_{23}^2 =
P_{12} P_{23} = P_2 $; that is, $ P_2 $ is also a projection, and $
P_1 + P_2 + P_3 = 1 $. It follows that $ P_2 H = \( (P_1+P_2) H \)
\cap \( (P_2+P_3) H \) = H_{12} \cap H_{23} = H_2 $. So, the relation
$ H = H_1 \oplus H_2 \oplus H_3 $ holds in the topological sense (that
is, when $ H $ is treated as a linear topological space). The
following lemma completes the proof.
\end{proof}

\begin{lemma}\label{3.2}
Let $ H $ be an FHS-space, $ H_1, H_2, H_3 $ its subspaces such that $
H = H_1 \oplus H_2 \oplus H_3 $ in the topological sense, and $
(H_1+H_2) \oplus H_3 = H = H_1 \oplus (H_2+H_3) $ in the FHS sense
(the bracketed sums being topological). Then $ H = H_1 \oplus H_2
\oplus H_3 $ in the FHS sense.
\end{lemma}

\begin{proof}
We introduce a Gaussian type space $ (\Om,\F,\P,G) $ and identify $ H
$ with $ G_0 $. Subspaces $ H_1, H_2, H_3 $ generate sub-\sif s $
\F_1, \F_2, \F_3 \subset \F $. Note that $ H_1+H_2 $ generates $ \F_1
\vee \F_2 $, the least \sif\ containing both $ \F_1 $ and $ \F_2
$. The following two lemmas (and one definition) complete the proof.
\end{proof}

\begin{definition}\label{3.3}
Let $ (\Om,\F,\P) $ be a measure type space, and $ \F_0, \F_1, \dots, \F_n
\subset \F $ sub-\sif s. We write $ \F_0 = \F_1 \otimes \dots \otimes
\F_n $, if $ \F_1 \vee \dots \vee \F_n = \F_0 $ and there exists $
P \in \P $ making $ \F_1, \dots, \F_n $ independent.\footnote{%
 Which means $ P ( A_1 \cap \dots \cap A_n ) = P (A_1) \dots P (A_n)
 $ for all $ A_1 \in \F_1, \dots, A_n \in \F_n $.}
\end{definition}

\begin{lemma}
Let $ (\Om,\F,\P,G) $ be a Gaussian type space, $ H = G_0 $ the
corresponding FHS-space, $ H_1, \dots, H_n \subset H $ subspaces, and
$ \F_1, \dots, \F_n \subset \F $ corresponding sub-\sif s. Then
\[
H = H_1 \oplus \dots \oplus H_n \quad \text{if and only if} \quad \F
= \F_1 \otimes \dots \otimes \F_n \, .
\]
\end{lemma}

\begin{proof}
`Only if': take a Gaussian measure $ \ga \in \Ga $ such that $
\|\cdot\|_\ga $ makes $ H_1, \dots, H_n $ orthogonal, then $ \ga $
makes $ \F_1, \dots, \F_n $ independent.

`If': some $ P \in \P $ makes $ \F_1, \dots, \F_n $ independent;
however, $ P $ need not be Gaussian. We take any Gaussian measure $
\ga_0 \in \Ga $ and introduce $ \F_k $-measurable densities
\[
f_k = \frac{ \ga_0 |_{\F_k} }{ P |_{\F_k} } \quad \text{for } k =
1,\dots,n \, .
\]
The measure $ \ga = f_1 \dots f_n \cdot P $ still makes $ \F_1,
\dots, \F_n $ independent, and $ \ga |_{\F_k} = f_k \cdot ( P
|_{\F_k} ) = \ga_0 |_{\F_k} $. So, w.r.t.\ $ \ga $ the spaces $ H_1,
\dots, H_n $ are \emph{independent} Gaussian spaces. Therefore their
sum $ H $ is Gaussian w.r.t.\ $ \ga $, that is, $ \ga \in \Ga $.
\end{proof}

\begin{lemma}\label{3.5}
Let $ (\Om,\F,\P) $ be a measure type space, and $ \F_1, \F_2, \F_3
\subset \F $ sub-\sif s such that
\[
(\F_1 \vee \F_2) \otimes \F_3 = \F = \F_1 \otimes (\F_2 \vee \F_3)
\, .
\]
Then
\[
\F = \F_1 \otimes \F_2 \otimes \F_3 \, .
\]
\end{lemma}

\begin{proof}
We have $ P,Q \in \P $ such that $ \F_1 \vee \F_2 $ and $ \F_3 $ are
$ P $-independent, while $ \F_1 $ and $ \F_2 \vee \F_3 $ are $ Q
$-independent. Consider the $ \F_1 \vee \F_2 $-measurable
density\footnote{%
 In fact, it is a conditional expectation w.r.t.\ $ P $, $ f_{12} =
 \cE{ \frac Q P }{ \F_1 \vee \F_2 } $.}
\[
f_{12} = \frac{ Q |_{\F_1\vee\F_2} }{ P |_{\F_1\vee\F_2} }
\]
and the measure $ R = f_{12} \cdot P \in \P $, then $ R
|_{\F_1\vee\F_2} = Q |_{\F_1\vee\F_2} $. The $ P $-independence of
$ \F_1 \vee \F_2 $ and $ \F_3 $ implies their $ R
$-independence. For every $ A \in \F_1 $, $ B \in \F_2 $, $ C \in \F_3
$
\begin{multline*}
R ( A \cap B \cap C ) = R ( A \cap B ) R (C) = Q ( A \cap B ) R (C) = \\
= Q (A) Q (B) R (C) = R (A) R (B) R (C) \, .
\end{multline*}
\end{proof}

The proof of Proposition \ref{3.1} is now complete. If you find the
proof of Lemma \ref{3.5} rather tricky, consider the following
calculation as a clue. Let $ \psi_k \in L_2 (\Om,\F_k,\P) $ ($ k =
1,2,3 $), then (writing for short $ \F_{12} = \F_1 \vee \F_2 $ etc.),
\begin{multline*}
(\psi_1\otimes\psi_2) \otimes \psi_3 = \frac{ \psi_1\otimes\psi_2 }{
  \sqrt{ P|_{\F_{12}} } } \cdot \frac{ \psi_3 }{ \sqrt{ P|_{\F_3} } }
  \cdot \sqrt P = \\
= \frac{ \psi_1 }{ \sqrt{ Q|_{\F_1} } } \cdot \frac{
  \psi_2 }{ \sqrt{ Q|_{\F_2} } } \cdot \sqrt{ \frac{ Q|_{\F_{12}} }{
  P|_{\F_{12}} } } \cdot \frac{ \psi_3 }{ \sqrt{ P|_{\F_3} } } \cdot
  \sqrt P \, ,
\end{multline*}
and the equality $ \| (\psi_1\otimes\psi_2) \otimes \psi_3 \|^2 = \|
\psi_1 \|^2 \| \psi_2 \|^2 \| \psi_3 \|^2 $ turns into the following
equality for functions $ f_1 = ( \psi_1 / \sqrt{Q|_{\F_1}} )^2 $,  $
f_2 = ( \psi_2 / \sqrt{Q|_{\F_2}} )^2 $,  $ f_3 = ( \psi_3 /
\sqrt{P|_{\F_3}} )^2 $:
\[
\int f_1 f_2 f_3 \frac{ Q|_{\F_{12}} }{ P|_{\F_{12}} } \, dP = \bigg(
\int f_1 \, dQ|_{\F_1} \bigg) \bigg( \int f_2 \, dQ|_{\F_2} \bigg)
\bigg( \int f_3 \, dP|_{\F_3} \bigg) \, .
\]
The argument is easily generalized for $ \F_1,\dots,\F_n $. Now we
turn to infinite sequences of subspaces.

\begin{definition}\label{3.6}
Let $ X $ be a metrizable topological space and $ X_1, X_2, \dotsc
\subset X $ closed subsets. We define $ \liminf_{n\to\infty} X_n $ as
the set of limits of all convergent sequences $ x_1, x_2, \dots $ such
that $ x_1 \in X_1, x_2 \in X_2, \dots $
\end{definition}

The set $ \liminf X_n $ is always closed.
If $ X $ is a linear topological space and $ X_n $
are linear subspaces, then $ \liminf X_n $ is a linear subspace. If $
X $ is the \sif\ $ \F $ of a measure type space $ (\Om,\F,\P) $ and $
X_n $ are sub-\sif s, then $ \liminf X_n $ is a sub-\sif.\footnote{%
 The natural topology of $ \F $ is defined by a metric $ \dist (A,B) =
 P ( A \setminus B ) + P ( B \setminus A ) $; the metric depends on $
 P \in \P $, but the topology does not. Of course, $ X $ is $ \F \bmod
 0 $ rather than $ \F $ itself.}
Proofs of these facts are left to the reader.

\begin{proposition}\label{3.7}
Let $ H $ be an FHS-space, $ E_n, F_n \subset H $ subspaces ($
n=1,2,\dots $), and $ \liminf E_n = H $. For each $ n $ denote by $
\cN_n $ the set of all admissible norms on $ H $ that make $ E_n, F_n
$ orthogonal. Then the set $ \liminf \cN_n $ is either the empty set,
or the whole $ \cN $.
\end{proposition}

\begin{proof}
We identify $ H $ with $ G_0 $ of a Gaussian type space $
(\Om,\F,\P,G) $ and use the natural homeomorphism $ \ga
\leftrightarrow \|\cdot\|_\ga $ between the space $ \cN $ of
admissible norms and the space $ \P_G $ of Gaussian
measures. Subspaces $ E_n, F_n $ generate corresponding sub-\sif s $
\E_n, \F_n \subset \F $. The set $ \cN_n \subset \cN $ corresponds to
$ \P_n \cap \P_G $, where $ \P_n $ consists of all measures making $
\E_n, \F_n $ independent. The following two lemmas complete the
proof.
\end{proof}

\begin{lemma}\label{3.8}
Let $ (\Om,\F,\P) $ be a measure type space, $ E, E_1, E_2, \dotsc
\subset L_0 (\P) $ closed linear subspaces, and $ \E, \E_1, \E_2,
\dotsc \subset \F $ the sub-\sif s generated by the subspaces. Then
\[
\liminf E_n = E \quad \text{implies} \quad \liminf \E_n \supset \E \,
.
\]
\end{lemma}

\begin{proof}
We have to prove that every $ f \in \liminf E_n $ is measurable
w.r.t.\ $ \liminf \E_n $. It suffices to prove that the set $ A = \{
\om : f(\om) \le a \} $ belongs to $ \liminf \E_n $ for every $ a $
such that the set $ \{ \om : f(\om) = a \} $ is negligible (indeed,
such $ a $ are dense in $ \R $). We take $ f_n \in E_n $ such that $
f_n \to f $ in $ L_0 (\P) $, consider sets $ A_n = \{ \om : f_n (\om)
\le a \} $ and note that $ A_n \in \E_n $ and $ A_n \to A $.
\end{proof}

Note. In general, $ \liminf \E_n $ need not be equal to $ \E $; it may
happen that $ \E_n = \F $ but $ \liminf E_n = \{ 0 \} $ (even for
one-dimensional $ E_n $). However, in the Gaussian case (that is, when
$ (\Om,\F,\P,G) $ is a Gaussian type space, and each $ E_n $ is a
Gaussian space) the equality $ \liminf \E_n = \E $ holds, which is
neither proved nor used here.

\begin{lemma}\label{3.9}
Let $ (\Om,\F,\P) $ be a measure type space, $ \E_n, \F_n \subset \F $
sub-\sif s ($ n = 1,2,\dots $), and $ \liminf \E_n = \F $. For each $
n $ denote by $ \P_n $ the set of all $ P \in \P $ such that $ \E_n $
and $ \F_n $ are $ P $-independent. Then

(a) the set $ \liminf \P_n $ is either the empty set or the whole $ \P
$;

(b) if $ \liminf \P_n = \P $ then $ \| (P-Q) |_{\F_n} \| \to 0 $ for
all $ P,Q \in \P $ (here $ \|\cdot\| $ means the total variation).
\end{lemma}

\begin{proof}
Assume that $ \liminf \P_n $ is nonempty; we have $ P_n, P \in \P $
such that $ P_n \to P $ and $ \E_n, \F_n $ are $ P_n
$-independent. Let $ K \subset (0,\infty) $ be a finite set and $ f :
\Om \to K $ an $ \F $-measurable function satisfying $ \int f \, dP =
1 $. Take $ \E_n $-measurable functions $ f_n : \Om \to K $ such that $
f_n \to f $ in $ L_0 (\P) $ and $ \int f_n \, dP = 1 $. The $ P_n
$-independent \sif s $ \E_n, \F_n $ are also $ Q_n $-independent,
where $ Q_n = f_n \cdot P_n \to f \cdot P = Q $ (since $ \| f_n \cdot
P_n - f_n \cdot P \| \le \| f_n \|_\infty \| P_n - P \| \to 0 $). Thus
$ Q_n \in \P_n $
and $ Q \in \liminf \P_n $. However, such measures $ Q $ (for all $ f
$ and $ K $) are dense in $ \P $. So, $ \liminf \P_n = \P $, which is
(a). Also, the $ P_n $-independence of $ \E_n, \F_n $ implies $ Q_n
|_{\F_n} = P_n |_{\F_n} $. However, $ \| (P_n-P) |_{\F_n} \| \le \|
P_n - P \| \to 0 $, and similarly $  \| (Q_n-Q) |_{\F_n} \| \to 0
$. So, $ \| (P-Q) |_{\F_n} \| \to 0 $ for all $ Q $ of a dense set,
therefore for all $ Q $, which is (b).
\end{proof}

The proof of Proposition \ref{3.7} is now complete.

\begin{definition}\label{3.10}
(a) Let $ H $ be an FHS-space, $ E_n, F_n \subset H $ subspaces, and $
\liminf E_n = H $. We say that $ F_n $ is asymptotically orthogonal to
$ E_n $, if there exists a convergent\footnote{%
 In the space $ \cN $ of all admissible norms, whose topology is
 defined in Sect.~2.}
sequence of admissible norms $ \|\cdot\|_n $ such that for each $ n $,
$ F_n $ is orthogonal to $ E_n $ w.r.t.\ $ \|\cdot\|_n $.

(b) Let $ (\Om,\F,\P) $ be a measure type space, $ \E_n, \F_n \subset
\F $ sub-\sif s, and $ \liminf \E_n = \F $. We say that $ \F_n $ is
asymptotically independent of $ \E_n $, if there exists a
convergent\footnote{%
 In the space $ \P $ whose topology is defined in Sect.~1.}
sequence of measures $ P_n \in \P $ such that for each $ n $, $ \F_n $
is independent of $ \E_n $ w.r.t.\ $ P_n $.
\end{definition}

Proposition \ref{3.7} states that the norms $ \|\cdot\|_n $ can be
chosen so as to converge to any given admissible norm, provided that $
F_n $ is asymptotically orthogonal to $ E_n $. Similarly, if $ \F_n $
is asymptotically independent of $ \E_n $, then the measures $ P_n $
can be chosen so as to converge to any given $ P \in \P $ due to
\ref{3.9}(a), and $ \| (P-Q) |_{\F_n} \| \to 0 $ for all $ P,Q \in \P
$ due to \ref{3.9}(b). Also, in the case of a Gaussian type space $
(\Om,\F,\P,G) $ and $ H = G_0 $, asymptotical orthogonality of
subspaces is equivalent to asymptotical independence of the
corresponding sub-\sif s.

\begin{lemma}\label{3.11}
Let $ H $ be an FHS-space, $ F_n \subset H $ subspaces, then the
following conditions are equivalent.

(a) For every $ f_1 \in F_1 $, $ f_2 \in F_2 $, \dots, if the sequence
$ (f_n) $ is bounded then $ f_n \to 0 $ weakly.\footnote{%
 That is, $ \ip{ f_n }{ h } \to 0 $ for every $ h \in H $.}

(b) For every finite-dimensional subspaces $ E_1, E_2, \dotsc $ such
that $ \liminf E_n = H $ there exist integers $ k_1 \le k_2 \le \dotsc
$ such that $ k_n \to \infty $ and $ F_n $ is asymptotically
orthogonal to $ E_{k_n} $.
\end{lemma}

\begin{proof}
We choose an admissible norm on $ H $, thus turning $ H $ into a
Hilbert space. Condition (a) becomes
\[
\forall h \in H \quad \sup_{f\in F_n, \|f\|\le1} \ip f h
\xrightarrow[n\to\infty]{} 0 \, ,
\]
that is,
\[
\forall h \in H \quad \angle (h,F_n) \xrightarrow[n\to\infty]{}
\frac\pi2 \, ,
\]
where the angle is defined by $ \cos \angle (h,F_n) = \sup \{ \ip h f
: f\in F_n, \|f\|\le1 \} $. It is equivalent to
\[
\forall E \quad \angle (E,F_n) \xrightarrow[n\to\infty]{} \frac\pi2 \, ,
\]
where $ E $ runs over finite-dimensional subspaces, and $ \cos \angle
(E,F_n) = \sup \{ \ip e f : e \in E, f \in F_n, \|e\| \le 1, \|f\| \le
1 \} $. The following lemma completes the proof, provided that $ k_n $
tends to $ \infty $ slowly enough. Namely, in terms of $ \de
(\cdot,\cdot) $ introduced there, it suffices that $ \de \( \angle
(E_{k_n},F_n), \dim E_{k_n} \) \xrightarrow[n\to\infty]{} 0 $.
\end{proof}

\begin{lemma}\label{3.12}
Let $ H $ be an FHS-space, $ E, F \subset H $ subspaces, $ \dim (E) <
\infty $, $ E \cap F = \{ 0 \} $. Then for every admissible norm $ \|
\cdot \|_1 $ there exists an admissible norm $ \|\cdot\|_2 $ such that
$ E,F $ are orthogonal w.r.t.\ $ \|\cdot\|_2 $, and
\[
\dist ( \|\cdot\|_1, \|\cdot\|_2 ) \le \de \( \angle(E,F), \dim E \)
\]
for some function $ \de : [0,\frac\pi2] \times \{0,1,2,\dots\} \to
(0,\infty) $ such that for every $ n $, $ \de(\al,n) \to 0 $ for $
\al \to \frac\pi2 $.\footnote{%
 Of course, $ \de $ does not depend on $ H,E,F $.}
\end{lemma}

\begin{proof}
We equip $ H $ with the norm $ \|\cdot\|_1 $, thus turning $ H $ into
a Hilbert space, and consider orthogonal projections $ Q_E, Q_F $ onto
$ E,F $ respectively. Introduce subspaces $ E \cap
F^\perp $, $ E^\perp \cap F $, $ E^\perp \cap F^\perp $ (here $
E^\perp $ is the orthogonal complement of $ E $); the subspaces are
orthogonal to each other, and invariant under both $ Q_E $ and $ Q_F
$. Therefore
\[
H = H_0 \oplus (E\cap F^\perp) \oplus (E^\perp\cap F)
\oplus (E^\perp\cap F^\perp) \, ,
\]
where $ H_0 $ is another subspace invariant under $ Q_E, Q_F $ (since
these operators are Hermitian). Introduce $ E_0 = E \cap H_0 $, $ F_0
= F \cap H_0 $, then $ Q_E h_0 = Q_{E_0} h_0 $ for all $ h_0 \in H_0 $
(since $ Q_E $ commutes with $ Q_{H_0} $), and $ Q_F h_0 = Q_{F_0} h_0
$. We may get rid of $ H_0^\perp $ by letting
\[
\| h_0 + h_1 \|_2^2 = \| h_0 \|_2^2 + \| h_1 \|_1^2 \quad \text{for
all $ h_0 \in H_0 $, $ h_1 \in H_0^\perp $.}
\]
In other words, we'll construct $ \|\cdot\|_2 $ on $ H_0 $ while
preserving both the given norm on $ H_1 $ and the orthogonality of $
H_0, H_1 $. Now we forget about $ H_1 $, assuming that $ H = H_0 $, $
E = E_0 $, $ F = F_0 $.

So, we have $ E \cap F = \{ 0 \} $, $ E \cap F^\perp = \{ 0 \} $, $
E^\perp \cap F = \{ 0 \} $, $ E^\perp \cap F^\perp = \{ 0 \} $. The
latter implies $ \dim (F^\perp) \le \codim (E^\perp) = \dim E
$. Similarly, $ \dim F \le \dim E $. Therefore $ H $ is
finite-dimensional, $ \dim H \le 2 \dim E $.

Both $ Q_E $ and $ Q_F $ commute with the Hermitian operator $ C =
\frac12 (2Q_E-1) (2Q_F-1) + \frac12 (2Q_F-1) (2Q_E-1) $. The spectrum
of $ C $ consists of some numbers $ \cos 2\phi_k $ of multiplicity 2
(though, some $ \phi_k $ may coincide), and $ 0 < \phi_k < \frac\pi2 $
(the case $ \phi_k = 0 $ is excluded by $ E \cap F = \{ 0 \} $; the
case $ \phi_k = \pi/2 $ is excluded by $  E \cap F^\perp = \{ 0 \} $,
$ E^\perp \cap F = \{ 0 \} $, $ E^\perp \cap F^\perp = \{ 0 \}
$). Accordingly, $ H_0 $ decomposes into the (orthogonal) direct sum
of planes, $ H = H_1 \oplus \dots \oplus H_d $, $ \dim H_k = 2 $,
invariant under $ Q_E, Q_F $. Subspaces $ E_k = E \cap H_k $, $ F_k =
F \cap H_k $ are two lines on the plane $ H_k $, and $ \angle
(E_k,F_k) = \phi_k $; $ k=1,\dots,d $; $ d \le \dim E $. Clearly,
\[
\angle (E,F) = \min ( \phi_1, \dots, \phi_d ) \, .
\]
We construct $ \|\cdot\|_2 $ on each $ H_k $ separately, while
preserving their orthogonality. Elementary 2-dimensional geometry
shows that the corresponding numbers $ \la'_k, \la''_k $ (two numbers
for each plane) are, in the optimal case,
\[
\la'_k = \( \tan \frac\phi2 \)^{-1/2} \, , \quad \la''_k = \( \tan
\frac\phi2 \)^{1/2} \, .
\]
The corresponding angle $ \be $ between $ \sqrt{\ga_1} $ and $
\sqrt{\ga_2} $ is given by
\[
\cos \be = \prod_{k=1}^d \bigg( \frac{ \tan^{-1/4} \frac\phi2 +
\tan^{1/4} \frac\phi2 }{ 2 } \bigg)^{-1} \, ,
\]
therefore
\[
\be \le \arccos \bigg( \frac{ \tan^{-1/4} \frac\al2 +
\tan^{1/4} \frac\al2 }{ 2 } \bigg)^{-\dim E} \, ;
\]
here $ \be = \angle (\sqrt{\ga_1}, \sqrt{\ga_2}) = \dist (
\|\cdot\|_1, \|\cdot\|_2 ) $, $ \al = \angle (E,F) $.
\end{proof}

\begin{definition}\label{3.13}
Let $ H $ be an FHS-space, $ F_n \subset H $ subspaces. We write $
\limsup F_n = \{0\} $, if the sequence $ (F_n) $ satisfies equivalent
conditions (a), (b) of Lemma \ref{3.11}.
\end{definition}

Note. More generally, one could define $ \limsup F_n $ as the set of
limits of all weakly convergent subsequences of all bounded sequences
$ f_1, f_2, \dotsc $ such that $ f_1 \in F_1, f_2 \in F_2, \dotsc $ It
is in general not a linear space, but anyway, $ ( \limsup E_n^\perp
)^\perp = \liminf E_n $ in a Hilbert space. For an FHS-space, as well
as a separable Banach space, $ F_n $ should be situated in the dual
space. However, all that is not needed here.

\begin{theorem}\label{3.14}
Let $ H $ be an FHS-space, $ E_n, F_n \subset H $ subspaces ($ n =
1,2,\dotsc $) such that $ \liminf E_n = H $, and $ \limsup F_n = \{0\}
$, and $ H = E_n \oplus F_n $ (in the FHS sense) for all $ n $. Then
there exist subspaces $ G_n, H_n \subset H $ such that
\[
E_n = G_n \oplus H_n \quad \text{and} \quad H = G_n \oplus H_n \oplus
F_n
\]
(both in the FHS sense), and $ \liminf G_n = H $, and $ H_n \oplus F_n
$ is asymptotically orthogonal to $ G_n $.
\end{theorem}

\begin{proof}
\begin{sloppypar}
We choose an admissible norm on $ H $, thus turning $ H $ into a
Hilbert space. Let $ L \subset H $ be a finite-dimensional subspace, $
L \ne \{0\} $. For any given $ n $ consider the pair $ L, E_n $. Its
geometry may be described (similarly to the proof of Lemma \ref{3.12})
via angles $ \phi_1^{(n)}, \dots, \phi_{d_n}^{(n)} \in [0,\frac\pi2)
$, $ d_n \le \dim L $. This time, zero angles are allowed, since $ L
\cap E_n $ need not be $ \{0\} $. It may happen that $ d_n < \dim L $,
since $ L \cap E_n^\perp $ need not be $ \{0\} $. However,
\[
\sup_{x\in L, x\ne0} \angle (x,E_n) = \al_n \to 0 \quad \text{for } n
\to \infty \, ;
\]
for large $ n $ we have $ \al_n < \pi/2 $ which implies $ d_n = d =
\dim L $ and $ \max ( \phi_1^{(n)}, \dots, \phi_d^{(n)} ) = \al_n
$. We may send $ L $ into $ E_n $ rotating it by $ \phi_1^{(n)},
\dots, \phi_d^{(n)} $. In other words, there is a rotation $ U_n : H
\to H $ such that
\[
U_n (L) \subset E_n \quad \text{and} \quad \| U_n - 1 \| \le 2 \sin
\frac{ \al_n }{ 2 } \xrightarrow[n\to\infty]{} 0 \, .
\]
\end{sloppypar}

We choose subspaces $ L_k \subset H $ such that $ \dim L_k = k $ and $
\liminf L_k = H $.\footnote{%
 Of course, one may take $ L_1 \subset L_2 \subset \dotsb $}
Introduce
\[
\al_{k,n} = \sup_{x\in L_k, x\ne0} \angle (x,E_n) \, ,
\]
then $ \al_{k,n} \xrightarrow[n\to\infty]{} 0 $ for each $ k $. On the
other hand, introduce
\[
\be_{k,n} = \frac\pi2 - \angle (L_k,F_n) \, .
\]
Similarly to the proof of Lemma \ref{3.11} we have $ \be_{k,n}
\xrightarrow[n\to\infty]{} 0 $ for each $ k $, therefore\footnote{%
 Recall that $ \de(\cdot,\cdot) $ is introduced in Lemma \ref{3.12}.}
$ \de ( \frac\pi2 - \be_{k_n,n}, k_n ) \xrightarrow[n\to\infty]{} 0 $
if $ k_n $ tends to $ \infty $ slowly enough. However, we choose $ k_1
\le k_2 \le \dotsb $, $ k_n \to \infty $ so as to satisfy a stronger
condition:
\[
\de \Big( \frac\pi2 - \al_{k_n,n} - \be_{k_n,n}, \, k_n \Big)
\xrightarrow[n\to\infty]{} 0 \, .
\]
We take
\[
G_n = U_n ( L_{k_n} ) \, ,
\]
where rotations $ U_n $ satisfy $ U_n (L_{k_n}) \subset E_n $ and $ \|
U_n - 1 \| \le 2 \sin (\frac12 \al_{{k_n},n} ) \to 0 $. Then $ \liminf
G_n = H $, and
\[
\frac\pi2 - \angle (G_n,F_n) \le \al_{k_n,n} + \be_{k_n,n} \, ;
\]
due to Lemma \ref{3.12}, $ F_n $ is asymptotically orthogonal to $ G_n
$. We take admissible norms $ \|\cdot\|_n \to \|\cdot\| $ such that $
F_n $ is orthogonal to $ G_n $ w.r.t.\ $ \|\cdot\|_n $. Consider the
orthogonal complement $ M_n $ of $ G_n $ w.r.t.\ $ \|\cdot\|_n $;
clearly, $ M_n $ is asymptotically orthogonal to $ G_n $. We have $
F_n \subset M_n $ and $ H = G_n \oplus M_n $ (in the FHS-sense). On
the other hand, $ G_n \subset E_n $ and $ H = E_n \oplus F_n
$. Proposition \ref{3.1} states that the subspace
\[
H_n = E_n \cap M_n
\]
satisfies $ G_n \oplus H_n \oplus F_n = H $ and $ G_n \oplus H_n = E_n
$ (and also $ H_n \oplus F_n = M_n $).
\end{proof}

%% file: sect4.tex
Recall the notion of a density matrix (borrowed from quantum
theory). Let $ H_1, H_2 $ be Hilbert spaces, $ H = H_1 \otimes H_2 $,
and $ \psi \in H $, $ \| \psi \| = 1 $. Every unit vector $ \xi \in
H_1 $ determines a subspace $ \xi \otimes H_2 \subset H $ and the
corresponding projection operator $ Q_{\xi\otimes H_2} = Q_\xi \otimes
\One_{H_2} $; here $ Q_\xi : H_1 \to H_1 $, $ Q_\xi x = (x,\xi) \xi $
is a one-dimensional projection, and $ \One_{H_2} : H_2 \to H_2 $, $
\One_{H_2} y = y $ the unit operator. The function
\[
\xi \mapsto \| Q_{\xi\otimes H_2} \psi \|^2
\]
is a quadratic form on $ H_1 $. The corresponding operator $ \rho_\psi
: H_1 \to H_1 $ satisfies
\[
\ip{ \rho_\psi \xi }{ \xi } = \| Q_{\xi\otimes H_2} \psi \|^2 = \ip{
( Q_\xi \otimes \One_{H_2} ) \psi }{ \psi }
\]
(that is, $ \langle \rho_\psi \rangle_\xi = \langle Q_\xi \otimes
\One_{H_2} \rangle_\psi $) for all $ \xi \in H_1 $; one calls $
\rho_\psi $ the \emph{density matrix} of $ \psi $ (on $ H_1 $). In
terms of an orthonormal basis $ (e_k) $ of $ H_2 $,
\begin{gather*}
\psi = \sum_k \psi_k \otimes e_k \quad \text{for some } \psi_k \in H_1
  \, ; \\
( Q_\xi \otimes \One_{H_2} ) \psi = \sum_k (Q_\xi \psi_k) \otimes e_k
  = \sum_k \ip{ \psi_k }{ \xi } \xi \otimes e_k \, ; \\
\ip{ \rho_\psi \xi }{ \xi } = \sum_k | \ip{ \psi_k }{ \xi } |^2 \, .
\end{gather*}
Note that
\begin{gather*}
\rho_\psi \ge 0 \, ; \qquad \Tr (\rho_\psi) = 1 \, ; \\
\Tr \( (\rho_{\psi_1}-\rho_{\psi_2}) A \) \le 2 \| \psi_1 - \psi_2 \|
  \cdot \| A \|
\end{gather*}
for all unit vectors $ \psi_1, \psi_2 \in H $ and operators $ A : H_1
\to H_1 $. The inequality may be proven as follows: $ \Tr ( \rho_\psi
Q_\xi ) = \Tr \( ( Q_\xi \otimes \One_{H_2} ) Q_\psi \) $; $ \Tr (
\rho_\psi A ) = \Tr \( ( A \otimes \One_{H_2} ) Q_\psi \) $; $ \Tr \(
( \rho_{\psi_1} - \rho_{\psi_2} ) A \) = \Tr \( ( A \otimes \One_{H_2}
) ( Q_{\psi_1} - Q_{\psi_2} ) \) \le \| A \otimes \One_{H_2} \| \cdot
\| Q_{\psi_1} - Q_{\psi_2} \| \le \| A \| \cdot 2 \| \psi_1 - \psi_2
\| $.

Note also that
\[
\psi' = ( U_1 \otimes U_2 ) \psi \quad \text{implies} \quad
\rho_{\psi'} = U_1 \rho_\psi U_1^*
\]
for all unitary operators $ U_1 : H_1 \to H_1 $, $ U_2 : H_2 \to H_2
$. Proof: $ \ip{ \rho_{\psi'} \xi }{ \xi } = \ip{ (Q_\xi\otimes
\One_{H_2}) \psi' }{ \psi' } = \ip{ (U_1 \otimes U_2)^* (Q_\xi\otimes
\One_{H_2}) (U_1\otimes U_2) \psi }{ \psi } = \ip{ ( U_1^* Q_\xi U_1
\otimes U_2^* \One_{H_2} U_2 ) \psi }{ \psi } = \ip{ ( Q_{U_1^* \xi}
\otimes \One_{H_2} ) \psi }{ \psi } = \ip{ \rho_\psi U_1^* \xi }{
U_1^* \xi } = \ip{ U_1 \rho_\psi U_1^* \xi }{ \xi } $ for all unit
vectors $ \xi \in H_1 $, $ \psi \in H $.

Assume in addition that $ H_1 = L_2 (\Om_1,\F_1,P_1) $, $ H_2 = L_2
(\Om_2,\F_2,P_2) $, then (after the usual identification) $ H = L_2
(\Om,\F,P) $ where $ (\Om,\F,P) = (\Om_1,\F_1,P_1) \otimes
(\Om_2,\F_2,P_2) $. We have
\[
\ip{ \rho_\psi \xi }{ \xi } = \iint \xi(\om') \overline{ \xi(\om'') }
\rho_\psi (\om',\om'') \, P_1 (d\om') P_1 (d\om'') \, ,
\]
where $ \rho_\psi $ is an element of $ L_2 \( (\Om_1,\F_1,P_1) \otimes
(\Om_1,\F_1,P_1) \) $ (the kernel of the operator $ \rho_\psi : H_1
\to H_1 $), namely,
\[
\rho_\psi ( \om', \om'' ) = \int \overline{ \psi(\om',\om_2) }
\psi(\om'',\om_2) \, P_2 (d\om_2) \, .
\]
In terms of a basis,
\begin{gather*}
\psi (\om_1,\om_2) = \sum_k \psi_k (\om_1) e_k (\om_2) \, ; \\
\rho_\psi (\om',\om'') = \sum_k \overline{ \psi_k(\om') }
\psi_k(\om'') \, .
\end{gather*}
The proof is basically a calculation:
\begin{multline*}
\rho_\psi (\om',\om'') = \int \overline{ \psi(\om',\om_2) }
  \psi(\om'',\om_2) \, P_2 (d\om_2) = \\
= \int \overline{ \sum_k \psi_k (\om') e_k (\om_2) } \sum_l \psi_l
  (\om'') e_l (\om_2) \, P_2 (d\om_2) = \\
= \sum_{k,l} \overline{ \psi_k (\om') } \psi_l (\om'') \int \overline{
  e_k (\om_2) } e_l (\om_2) \, P_2 (d\om_2) = \sum_k \overline{ \psi_k
  (\om') } \psi_k (\om'') \, ;
\end{multline*}
\begin{multline*}
\ip{ \rho_\psi \xi }{ \xi } = \sum_k | \ip{ \psi_k }{ \xi } |^2 =
  \sum_k \Big| \int \psi_k (\om) \overline{ \xi(\om) } \, P_1(d\om)
  \Big|^2 = \\
= \sum_k \Big( \int \overline{ \psi_k(\om') } \xi (\om') \, P_1
  (d\om') \cdot \int \psi_k (\om'') \overline{ \xi(\om'') } \, P_1
  (d\om'') \Big) = \\
= \iint \xi(\om') \overline{ \xi(\om'') } \Big( \sum_k \overline{
  \psi_k(\om') } \psi_k(\om'') \Big) \, P_1(d\om') P_1(d\om'') = \\
= \iint \xi (\om') \overline{ \xi (\om'') } \rho_\psi (\om',\om'') \,
  P_1(d\om') P_1(d\om'') \, .
\end{multline*}

We turn to measure type spaces $ (\Om,\F,\P) = (\Om_1,\F_1,\P_1)
\otimes (\Om_2,\F_2,\P_2) $ (which means $ P_1 \otimes P_2 \in \P $
for some, therefore all, $ P_1 \in \P_1 $ and $ P_2 \in \P_2 $), and
the corresponding Hilbert spaces; $ H = L_2 (\Om,\F,\P) = L_2
(\Om_1,\F_1,\P_1) \otimes L_2 (\Om_2,\F_2,\P_2) = H_1 \otimes H_2 $
under the natural identification
\[
\frac{ \psi_1 \otimes \psi_2 }{ \sqrt{ P_1 \otimes P_2 } }
(\om_1,\om_2) = \frac{ \psi_1 }{ \sqrt{P_1} } (\om_1) \cdot \frac{
\psi_2 }{ \sqrt{P_2} } (\om_2) \, ,
\]
that is,
\[
\psi_1 \otimes \psi_2 = \bigg( \frac{ \psi_1 }{ \sqrt{P_1} } \otimes
\frac{ \psi_2 }{ \sqrt{P_2} } \bigg) \cdot \sqrt{ P_1 \otimes P_2 }
\in L_2 (\P)
\]
for all $ \psi_1 \in L_2 (\P_1) $, $ \psi_2 \in L_2 (\P_2) $, $ P_1
\in \P_1 $, $ P_2 \in \P_2 $. An element $ \psi \in H $, $ \| \psi \|
= 1 $, determines an operator $ \rho_\psi : H_1 \to H_1 $, whose
kernel (denoted also by $ \rho_\psi $) belongs to $ L_2 \(
(\Om_1,\F_1,\P_1) \otimes (\Om_1,\F_1,\P_1) \) $;
\begin{gather*}
\ip{ \rho_\psi \xi }{ \xi } = \iint \frac{ \xi }{ \sqrt{P_1} } (\om')
  \frac{ \xi }{ \sqrt{P_1} } (\om'') \frac{ \rho_\psi }{ \sqrt { P_1
  \otimes P_1 } } (\om', \om'') \, P_1 (d\om') P_1 (d\om'') \, , \\
\frac{ \rho_\psi }{ \sqrt { P_1 \otimes P_1 } } (\om', \om'') = \int
  \overline{ \frac{ \psi }{ \sqrt{P_1\otimes P_2} } (\om',\om_2) }
  \cdot \frac{ \psi }{ \sqrt{P_1\otimes P_2} } (\om'',\om_2) \, P_2
  (d\om_2)
\end{gather*}
for all $ \xi \in L_2 (\P_1) $, $ \| \xi \| = 1 $, and $ P_1 \in \P_1
$, $ P_2 \in \P_2 $.

Let $ (\Om,\F,\P) $ be a measure type space and $ \F_1,\F_2 \subset \F
$ sub-\sif s such that $ \F = \F_1 \otimes \F_2 $ (as defined by
\ref{3.3}), then
\[
L_2 (\F) = L_2 (\F_1) \otimes L_2 (\F_2) \, .
\]
Here $ L_2 (\F) = L_2 (\Om,\F,\P) $ and $ L_2 (\F_1) = L_2
(\Om,\F_1,\P|_{\F_1}) $, the same for $ \F_2 $,\footnote{%
 Of course, $ \P|_{\F_1} = \{ P|_{\F_1} : P \in \P \} $ consists of
 restricted measures. Alternatively one may introduce quotient spaces
 $ \Om_1 = \Om / \F_1 $, $ \Om_2 = \Om / \F_2 $ and identify $ \Om $
 with $ \Om_1 \times \Om_2 $.}
and the natural identification is made, namely,
\[
\( f_1 \cdot \sqrt{ P |_{\F_1} } \) \otimes \( f_2 \cdot \sqrt{ P
|_{\F_2} } \) = (f_1 f_2) \cdot \sqrt{ P }
\]
whenever $ f_1 \in L_2 (\Om,\F_1,P) $, $ f_2 \in L_2 (\Om,\F_2,P) $,
and $ P \in \P $ makes $ \F_1, \F_2 $ independent. We need a
counterpart of Lemma \ref{3.9}.

\begin{lemma}\label{4.1}
Let $ (\Om,\F,\P) $ be a measure type space, $ \E_n,\F_n \subset \F $
sub-\sif s ($ n = 1,2,\dots $), $ \F = \E_n \otimes \F_n $ for each $
n $, and $ \liminf \E_n = \F $, and $ \F_n $ is asymptotically
independent of $ \E_n $. Then for every $ P \in \P $
\[
\liminf \( L_2 (\E_n) \otimes \sqrt{ P |_{\F_n} } \) = L_2 (\F) \, .
\]
That is, for every $ \psi \in L_2 (\F) $ there exist $ \xi_n \in L_2
(\E_n) $ such that $ \| \psi - \xi_n \otimes \sqrt{ P |_{\F_n} } \|
\to 0 $ when $ n \to \infty $.
\end{lemma}

\begin{proof}
Similarly to the proof of \ref{3.9}, we take $ P_n \in \P $ such that
$ \E_n, \F_n $ are $ P_n $-independent and $ P_n \to P $. We consider
an arbitrary finite set $ K \subset \R $, an arbitrary $ \F
$-measurable function $ f : \Om \to K $, and the corresponding vector
$ \psi = f \sqrt P \in L_2 (\F) $. We construct $ \E_n $-measurable
functions $ f_n : \Om \to K $ such that $ f_n \to f $ in $ L_0 (\P) $,
and corresponding vectors $ \psi_n = f_n \cdot \sqrt{ P_n } \in L_2
(\F) $. Independence of $ \E_n, \F_n $ w.r.t.\ $ P_n $ means that $
\sqrt{ P_n } = \sqrt{ P_n |_{\E_n} } \otimes \sqrt{ P_n |_{\F_n} } $,
therefore $ \psi_n = \xi_n \otimes \sqrt{ P_n |_{\F_n} } \in L_2
(\E_n) \otimes \sqrt{ P_n |_{\F_n} } $, where $ \xi_n = f_n \cdot
\sqrt{ P_n |_{\E_n} } \in L_2 (\E_n) $. However, $ \psi_n \to \psi $,
since $ \| f_n \sqrt{ P_n } - f_n \sqrt P \| \le \| f_n \|_\infty \|
\sqrt{ P_n } - \sqrt P \| \to 0 $. Also $ \| P_n |_{\F_n} - P |_{\F_n}
\| \le \| P_n - P \| \to 0 $, therefore $ \| \sqrt{ P_n |_{\F_n} } -
\sqrt{ P |_{\F_n} } \| \to 0 $. So, $ \| \psi - \xi_n \otimes \sqrt{ P
|_{\F_n} } \| \le \| \psi - \psi_n \| + \| \xi_n \otimes \sqrt{ P_n
 |_{\F_n} } - \xi_n \otimes \sqrt{ P |_{\F_n} } \| \to 0 $, and $ \psi
\in \liminf \( L_2 (\E_n) \otimes \sqrt{ P |_{\F_n} } \) $. It remains
to note that such vectors $ \psi $ (for all $ f $ and $ K $) are dense
in $ L_2 (\F) $.
\end{proof}

%% file: sect5.tex
Recall the correspondence (described in Sect.\ 2) between Gaussian
type spaces $ (\Om,\F,\P,G) $ and FHS-spaces $ G_0 = G / \Const $. We
know that any FHS-space $ G_0 $ determines (up to isomorphism) the
corresponding Gaussian type space $ (\Om,\F,\P,G) $, which in turn
determines the Hilbert space $ H = L_2 (\Om,\F,\P) $. I denote the
relation by
\[
H = \Exp (G_0) \, ,
\]
and call $ H $ the Fock exponential of $ G_0 $.\footnote{%
 By choosing an admissible norm on $ G_0 $ one turns $ G_0 $ into a
 Hilbert space, in which case $ \Exp (G_0) $ contains a special
 element, `ground state vector' $ \sqrt \ga $ (where $ \ga $ is the
 Gaussian measure corresponding to the chosen norm); that is the
 classical Fock construction.}
Why call it `exponential'? Since
\[
G_0 = G_1 \oplus G_2 \quad \text{implies} \quad \Exp (G_0) = \Exp
(G_1) \otimes \Exp (G_2)
\]
in the following sense. Assume that $ G_1,G_2 \subset G_0 $ are
subspaces such that $ G_0 = G_1 \oplus G_2 $ (in the FHS sense, as
defined in Sect.~3). Then the sub-\sif s $ \F_1,\F_2 \subset \F $,
generated by $ G_1,G_2 $ respectively, satisfy $ \F = \F_1 \otimes
\F_2 $, therefore $ L_2 (\F) = L_2 (\F_1) \otimes L_2 (\F_2) $, as
explained in Sect.~4. So, every decomposition of $ G_0 $ into a direct
sum determines a decomposition of $ \Exp (G_0) $ into a tensor
product.

Given such a decomposition $ G_0 = G_1 \oplus G_2 $, every unit vector
$ \psi \in \Exp (G_0) $ determines a density matrix $ \rho_\psi $ on $
\Exp (G_2) $. Do not think, however, that $ \rho_\psi $ is uniquely
determined by $ \psi $ and $ G_2 $; also $ G_1 $ influences $
\rho_\psi $.

Consider an infinite sequence of decompositions, $ G_0 = E_n \oplus
F_n $, of a single FHS-space $ G_0 $. We want to know, whether or not
$ \| \rho_n (\psi_1) - \rho_n (\psi_2) \| \to 0 $ when $ n \to \infty
$ for all unit vectors $ \psi_1, \psi_2 \in \Exp (G_0) $; here $
\rho_n (\psi) $ is the density matrix on $ \Exp (F_n) $ that
corresponds to $ \psi \in \Exp (E_n) \otimes \Exp (F_n) $, and the
trace norm is used,
\[
\| \rho_n (\psi_1) - \rho_n (\psi_2) \| = \sup_{\|A\|\le1} \Tr \(
\rho_n (\psi_1) - \rho_n (\psi_2) ) A \) \, ;
\]
here $ A $ runs over Hermitian operators on $ \Exp (F_n) $. In general
we have only $ \| \rho_n (\psi_1) - \rho_n (\psi_2) \| \le 2 \| \psi_1
- \psi_2 \| $.

Assume that $ \liminf E_n = G_0 $. If $ F_n $ is asymptotically
orthogonal to $ E_n $, then $ \| \rho_n (\psi_1) - \rho_n (\psi_2) \|
\to 0 $, which follows easily from Lemma \ref{4.1}. Namely, the lemma
represents $ \psi $ as the limit of $ \xi_n \otimes \sqrt{ P |_{\F_n}
} $; however, $ \rho_n ( \xi_n \otimes \sqrt{ P |_{\F_n} } ) $ is the
one-dimensional projection $ Q_n $ onto $ \sqrt{ P |_{\F_n} } $,
therefore $ \| \rho_n (\psi) - Q_n \| \le 2 \| \psi - \xi_n \otimes
\sqrt{ P |_{\F_n} } \| \to 0 $, and so, $ \| \rho_n (\psi_1) - \rho_n
(\psi_2) \| \le \| \rho_n (\psi_1) - Q_n \| + \| \rho_n (\psi_2) - Q_n
\| \to 0 $. However, the asymptotical orthogonality condition may be
dropped, as we'll see now.

\begin{proposition}\label{5.1}
Let $ G_0 $ be an FHS-space, $ E_n, F_n \subset G_0 $ subspaces such
that $ \liminf E_n = G_0 $, and $ \limsup F_n = \{0\} $, and $ G_0 =
E_n \oplus F_n $ (in the FHS sense) for all $ n $. Then $ \| \rho_n
(\psi_1) - \rho_n (\psi_2) \| \to 0 $ when $ n \to \infty $, for all
unit vectors $ \psi_1, \psi_2 \in \Exp (G_0) $.
\end{proposition}

\begin{proof}
Theorem \ref{3.14} gives us $ G_n, H_n \subset G_0 $ such that $ E_n =
G_n \oplus H_n $ and $ G_0 = G_n \oplus H_n \oplus F_n $ (both in the
FHS sense), and $ \liminf G_n = G_0 $, and $ H_n \oplus F_n $ is
asymptotically orthogonal to $ G_n $. Lemma \ref{4.1} gives us a
representation
\[
\psi = \lim_{n\to\infty} ( \xi_n \otimes \chi_n )
\]
for an arbitrary unit vector $ \psi \in \Exp (G_0) $; here $ \xi_n $
are unit vectors of $ \Exp (G_n) $, $ \chi_n $ are unit vectors of $
\Exp (H_n \oplus F_n) $, and these $ \chi_n $ (unlike $ \xi_n $) do
not depend on $ \psi $. We have $ \| \rho_n (\psi) - \rho_n (\xi_n
\otimes \chi_n) \| \le 2 \| \psi - \xi_n \otimes \chi_n \| \to 0
$. However, the density matrix $ \rho_n (\xi_n \otimes \chi_n) $ on $
\Exp (F_n) $ that corresponds to the vector $ \xi_n \otimes \chi_n \in
\Exp (G_n) \otimes \Exp (H_n) \otimes \Exp (F_n) $ is the same as the
density matrix on $ \Exp (F_n) $ that corresponds to the vector $
\chi_n \in \Exp (H_n) \otimes \Exp (F_n) $. The vector does not depend
on $ \psi $, therefore $ \rho_n (\xi_n \otimes \chi_n) $ does not
depend on $ \psi $. So, $ \| \rho_n (\psi_1) - \rho_n (\psi_2) \| \le
\| \rho_n (\psi_1) - \rho_n (\xi_n \otimes \chi_n) \| + \| \rho_n
(\psi_2) - \rho_n (\xi_n \otimes \chi_n) \| \to 0 $.
\end{proof}

Consider a decomposition $ G_0 = G_1 \oplus G_2 $ of an FHS-space $
G_0 $ into the direct sum (in the FHS sense) of subspaces $ G_1, G_2
\subset G_0 $, and the corresponding decomposition $ G^0 = G^1 \oplus
G^2 $ of its dual FHS-space $ G^0 $; that is, $ G^1 = G_2^\perp $ is
the annihilator of $ G_2 $ in $ G^0 $, and $ G^2 = G_1^\perp $. Of
course, also $ G_1 = (G^2)^\perp $ and $ G_2 = (G^1)^\perp
$. Introduce an admissible norm $ \|\cdot\| $ on $ G_0 $ (note that $
G_1, G_2 $ need not be orthogonal w.r.t.\ $ \|\cdot\| $), and its dual
norm on $ G^0 $ (denoted by $ \|\cdot\| $ as well). Consider
\begin{align*}
\dist ( g, G_2 ) &= \inf_{g_2\in G_2} \| g - g_2 \| = \sup_{x_1\in
  G^1, \|x_1\|\le1} \ip{ g }{ x_1 } \, , \\
\dist ( x, G^2 ) &= \inf_{x_2\in G^2} \| x - x_2 \| = \sup_{g_1\in
  G_1, \|g_1\|\le1} \ip{ g_1 }{ x }
\end{align*}
for any $ g \in G_0 $, $ x \in G^0 $.

Introduce the corresponding Hilbert spaces $ H = \Exp (G_0) $, $ H_1 =
\Exp (G_1) $, $ H_2 = \Exp (G_2) $; we have $ H = H_1 \otimes H_2
$. Recall the operators $ U_x $ and $ V_g $ for $ x \in G^0 $, $ g \in
G_0 $, satisfying the Canonical Commutation Relations $ V_g U_x =
e^{i\ip g x} U_x V_g $. Let $ (\Om,\F,\P,G) $ be the corresponding
Gaussian type space (thus, $ H = L_2(\Om,\F,\P) $), and $ \ga \in \P_G
$ a Gaussian measure such that $ \|\cdot\|_\ga = \|\cdot\| $. We know
that the vector $ \psi = \sqrt\ga \in H $ satisfies $ \ip{ U_x \psi }{
\psi } = \exp \( -\frac18 \|x\|^2 \) $, $ \ip{ V_g \psi }{ \psi } =
\exp \( -\frac12 \|g\|^2 \) $ for all $ x \in G^0 $, $ g \in G_0
$. Denote by $ \rho (\psi) $ the corresponding density matrix on $ H_1
$. In the following lemma, $ \rho(\psi') $ for some other $ \psi' $
are the corresponding density matrices on $ H_1 $, and a function $ M :
[0,\infty) \to [0,\infty) $ is defined by\footnote{%
 The exact form of the function is of no importance here; we only need
 to know that $ M(r_n) \to 0 $ implies $ r_n \to 0 $.}
\[
M(r) = \max_{\phi\in[0,\pi]} \bigg( \exp \Big( -\frac{\phi^2}{2r^2}
\Big) \cdot 2 \sin \frac\phi2 \bigg) \, .
\]

\begin{lemma}\label{5.2}
For all $ x \in G^0 $, $ g \in G_0 $
\begin{align*}
\| \rho (\psi) - \rho (U_x \psi) &\| \ge M \( \dist (x,G^2) \) \, , \\
\| \rho (\psi) - \rho (V_g \psi) &\| \ge M \( 2 \dist (g,G_2) \) \,
 . \\
\end{align*}
\end{lemma}

\begin{proof}
Let $ x = y + z $, $ y \in G^1 $, $ z \in G^2 $, then $ U_x =
U_y^{(1)} \otimes U_z^{(2)} $ (the notation being self-explanatory),
which implies $ \rho (U_x \psi) = U_y^{(1)} \rho(\psi) U_{-y}^{(1)}
$. For every $ g \in G_1 $
\[
\Tr \( \rho(U_x\psi) V_g^{(1)} \) = \Tr \( U_y^{(1)} \rho(\psi)
U_{-y}^{(1)} V_g^{(1)} \) = \Tr \( \rho(\psi) U_{-y}^{(1)} V_g^{(1)}
U_y^{(1)} \) \, ;
\]
however, $ V_g^{(1)} U_y^{(1)} = \exp (i\ip g y) U_y^{(1)} V_g^{(1)} $
and $ \ip g y = \ip g x $ (since $ g \in (G^2)^\perp $); we have $ 
\Tr \( \rho(U_x\psi) V_g^{(1)} \) = \exp (i\ip g x) \Tr \(
\rho(\psi) V_g^{(1)} \) $. Note that $ | \Tr \( ( \rho(U_x\psi) -
\rho(\psi) ) V_g^{(1)} \) | \le \| \rho(U_x\psi) - \rho(\psi) \| $,
since $ \Re \( e^{i\al} \Tr \( ( \rho(U_x\psi) - \rho(\psi) )
V_g^{(1)} \) = \Tr \( ( \rho(U_x\psi) - \rho(\psi) ) \Re ( e^{i\al}
V_g^{(1)} ) \) \le \| \rho(U_x\psi) - \rho(\psi) \| $ for all $ \al
$. We have
\begin{gather*}
\| \rho(U_x\psi) - \rho(\psi) \| \ge | 1 - e^{i\ip g x} | \cdot | \Tr
  \( \rho(\psi) V_g^{(1)} \) | \, ; \\
| 1 - e^{i\ip g x} | = 2 \Big| \sin \frac{\ip g x}2 \Big| \, ; \\
\Tr \( \rho(\psi) V_g^{(1)} \) = \Tr \( ( V_g^{(1)} \otimes \One_{H_2}
  ) Q_\psi \) = \ip{ ( V_g^{(1)} \otimes \One_{H_2} ) \psi }{ \psi } =
  \qquad \qquad \\
\qquad \qquad = \ip{ V_g \psi }{ \psi } = \exp \( -\tfrac12 \| g \|^2
  \) \, ;
\end{gather*}
so,
\[
\| \rho(U_x\psi) - \rho(\psi) \| \ge \sup_{g\in G_1} \bigg( \exp \(
-\tfrac12 \| g \|^2 \) \cdot 2 \Big| \sin \frac{ \ip g x }{ 2 } \Big|
 \bigg) \, .
\]
Denote $ r = \dist (x,G^2) $. We need only one ray of vectors $ g \in
G_1 $ such that $ \ip g x = r \|g\| $. For every $ \phi \in [0,\pi] $
there exists such $ g $, satisfying $ \ip g x = \phi $ and $ \| g \| =
\phi / r $, which gives $ \| \rho(U_x\psi) - \rho(\psi) \| \ge \exp \(
-\frac{ \phi^2 }{ 2r^2 } \) \cdot 2 \sin \frac\phi2 $; the supremum
over $ \phi $ gives the first inequality.

For the second inequality, the proof is quite similar. Only $
U_y^{(1)} $ and $ V_g^{(1)} $ change places, and $ \exp \( -\frac18
\|x\|^2 \) $ appears instead of $ \exp \( -\frac12 \|g\|^2 \) $, which
leads to $ M(2r) $ instead of $ M(r) $.
\end{proof}

\begin{theorem}\label{5.3}
Let $ G_0 $ be an FHS-space, $ E_n, F_n \subset G_0 $ subspaces such
that $ G_0 = E_n \oplus F_n $ (in the FHS sense) for all $ n $. For
every unit vector $ \psi \in \Exp (G_0) $ let $ \rho_n (\psi) $ denote
the density matrix on $ F_n $ that corresponds to $ \psi $. Then the
following two conditions are equivalent.

(a) $ \liminf E_n = G_0 $ and $ \limsup F_n = \{ 0 \} $.

(b) $ \| \rho_n (\psi_1) - \rho_n (\psi_2) \| \to 0 $ when $ n \to
\infty $, for all unit vectors $ \psi_1, \psi_2 \in \Exp (G_0) $.
\end{theorem}

\begin{proof}
Proposition \ref{5.1} gives (a) \imp (b). Assume (b); we have to prove
(a). We choose a Gaussian measure $ \ga $ and apply Lemma \ref{5.2} to
$ \psi = \sqrt\ga $:
\begin{align*}
M \( \dist ( x, F_n^\perp ) \) \le \| \rho_n (\psi) - \rho_n (U_x
  \psi) \| \to 0 \, , \\
M \( 2 \dist ( g, E_n ) \) \le \| \rho_n (\psi) - \rho_n (V_g \psi) \|
  \to 0 \, ,
\end{align*}
which implies that $ \dist ( x, F_n^\perp ) \to 0 $ and $ \dist ( g,
E_n ) \to 0 $ (when $ n \to \infty $) for all $ g \in G_0 $, $ x \in
G^0 $ (the dual to $ G_0 $). The latter, $\displaystyle \inf_{e\in
E_n} \| g - e \| \to 0 $, shows that $ \liminf E_n = G_0 $. The
former, $\displaystyle \sup_{f\in F_n, \|f\|\le1} \ip f x \to 0 $,
shows that $ \limsup F_n = \{ 0 \} $.
\end{proof}

%% file: sect6.tex
Recall some notions and results about Borel measurability (see
\cite[Chapter 3]{Ar76}, \cite[Chapter 3]{Sr}). A \emph{Borel space} is
a set equipped with a \sif\ (of subsets). The subsets belonging to the
\sif\ are called Borel measurable sets (or `measurable sets', or
`Borel sets'). A subset of a Borel space is naturally a Borel
space. The product of two Borel spaces is naturally a Borel space. A
Borel measurable map (or `measurable map', or `Borel map')
is a map from one Borel space into another, such that the inverse
image of every measurable set is a measurable set. A Borel isomorphism
between two Borel spaces is an invertible measurable map whose inverse
is also measurable. (Note that no measure (type) is given, and so, no
subset is negligible; every single point counts.)

A \emph{Polish space} is a topological space which is homeomorphic to
a separable complete metric space. A Polish space is naturally
equipped with the \sif\ generated by all open sets, thus, it is a
Borel space. Surprisingly, the Borel space does not depend (up to
isomorphism) on the Polish space, as far as it is uncountable. That
is, every two uncountable Polish spaces are Borel
isomorphic. Moreover, all uncountable Borel sets in Polish spaces are
Borel isomorphic. A Borel space is called \emph{standard,} if it is
isomorphic to a Borel subset of a Polish space. Up to isomorphism,
there is a single uncountable standard Borel space, a single countable
(infinite) one, and for each (finite) $ n $, a single $ n $-point one.

Let $ X $ be a Polish space and $ \bF (X) $ the set of all nonempty
closed subsets of $ X $. There is a natural Borel structure on $ \bF
(X) $, namely, the \sif\ generated by sets of the form
\[
\{ F \in \bF (X) : F \cap U \neq \emptyset \}
\]
where $ U $ varies over open sets in $ X $. Thus $ \bF (X) $ is a
Borel space; it is called the Effros Borel space of $ X $, and is
standard (see \cite[Th.~3.3.10]{Sr}). There is a sequence $ (f_n) $ of
Borel measurable maps $ f_n : \bF (X) \to X $ such that every $ F \in
\bF(X) $ is the closure of the countable (or finite) set $ \{ f_1(F),
f_2(F), \dots \} $; it is called Castaing's theorem (see
\cite[Prop.~5.2.7]{Sr}). Therefore, for any Borel space $ T $, a
general form of a measurable map $ f : T \to \bF (X) $ is
\begin{equation}\label{6.1}
f(t) = \closure \( \{ f_1(t), f_2(t), \dots \} \) \, ,
\end{equation}
where $ f_1, f_2, \dotsc : T \to X $ are measurable maps.
Note that the disjoint union of all $ F \in \bF (X) $, defined as the
set of all pairs $ (F,x) $ such that $ F \in \bF (X) $ and $ x \in F
$, is naturally a standard Borel space, since it is a Borel subset of
$ \bF (X) \times X $.

Now we apply all that to our matter. Let $ H $ be a (separable)
Hilbert space and $ \bL (H) $ the set of all (closed) linear subspaces
of $ H $. Then $ H $ is a Polish space, $ \bL (H) \subset \bF (H) $,
and we get measurable maps $
f_n : \bL (H) \to H $ such that every $ L \in \bL (H) $ is spanned by
(and even the closure of) $ \{ f_1(L), f_2(L), \dots \} $. Applying
the usual orthogonalization process to the sequence $ \( f_n (L) \) $
we get a new sequence, denote it again by $ \( f_n(L) \) $, such that
\begin{equation}\label{6.2}
\{ f_n(L) : 1 \le n < 1+\di L \} \quad \text{is an orthonormal basis
of } L
\end{equation}
for every $ L \in \bL (H) $, and still, $ f_n : \bL (H) \to H $ are
Borel measurable, and in addition, $ f_n (L) = 0 $ when
$ n \ge 1+\di L $. The same argument shows also that $ \di L \in \{
0,1,2,\dots; \infty \} $ is a
measurable function of $ L $.

There is a natural map, $ F \mapsto \spa F $, from $ \bF (H) $ to $
\bL (H) $; namely,
$ \spa F $ is the (closed) subspace spanned by $ F
$. The map is measurable, which follows from \eqref{6.1}. Indeed, if $
F $ is the closure of $ \{ f_n (F) : n = 1,2,\dots \} $, then $ \spa
F $ is the closure of
\begin{equation}\label{6.3}
\{ \al_1 f_1 (F) + \dots + \al_n f_n (F) : \al_1,\dots,\al_n \in \Q,
\, n = 1,2,\dots \} \, ,
\end{equation}
still a countable set of measurable functions (indexed by finite
sequences $ (\al_1,\dots,\al_n) $ of rational numbers).

\begin{lemma}\label{6.4}
Linear subspaces of a Hilbert space are a standard Borel
space.\footnote{%
 \emph{Closed} linear subspaces of a \emph{separable} Hilbert space,
 of course.}
\end{lemma}

\begin{proof}
$ \bL (H) = \{ F \in \bF (H) : \spa (F) = F \} $, therefore $ \bL (H)
$ is a Borel subset of the standard Borel space $ \bF (H)
$.\footnote{%
 Indeed, $ F \mapsto ( F, \spa F ) $ is a Borel map $ \bF (H) \to \bF
 (H) \times \bF (H) $, and the diagonal is a Borel subset of $ \bF (H)
 \times \bF (H) $.}
\end{proof}

\begin{lemma}\label{6.5}
The closure of $ L_1 + L_2 $ is a jointly Borel measurable function of
linear subspaces $ L_1, L_2 $ of a Hilbert space.
\end{lemma}

The proof is left to the reader. Hint: Similar to \eqref{6.3} but
simpler.\footnote{%
 An alternative way: $ L_1 \cup L_2 $ is measurable in $ L_1, L_2 $ by
 \cite[Exercise 3.3.11(ii)]{Sr}, thus $ \spa (L_1 \cup L_2) $ is also
 measurable. It is a good luck that we need unions, not intersections;
see the note after Exercise 3.3.11 in \cite{Sr}.}

\begin{lemma}\label{6.6}
The orthogonal projection of $ x $ to $ L $ is jointly measurable in $
x \in H $ and $ L \in \bL (H) $.
\end{lemma}

\begin{proof}
The projection is the limit (for $ n \to \infty $) of the orthogonal
projection of $ x $ to $ \spa \{ f_k (L) : k \le n \}
$. Measurability of the latter implies that of the former.
\end{proof}

Note also that the disjoint union of all $ L \in \bL (H) $ is a
standard Borel space, and linear operations are measurable in the
following sense: $ (L,h_1+h_2) $ is jointly Borel measurable in $
(L,h_1) $ and $ (L,h_2) $ on the domain consisting of all pairs $ \(
(L_1,h_1), (L_2,h_2) \) $ where $ L_1,L_2 \in \bL (X) $, $ h_1 \in L_1
$, $ h_2 \in L_2 $ satisfy $ L_1 = L_2 $. The same for other linear
combinations, and for the scalar product.

Let $ H $ be an FHS-space (rather than a Hilbert space). Lemmas
\ref{6.4} and \ref{6.5} still hold.

\begin{lemma}
The set of all pairs $ (L_1,L_2) $ such that $ H = L_1 \oplus L_2 $
(in the FHS sense, see Sect.~3) is a Borel subset of $ \bL (H) \times
\bL (H) $.
\end{lemma}

\begin{proof}
The relation $ H = L_1 \oplus L_2 $ means that, first, $ L_1, L_2 $
are orthogonal in some admissible norm, and second, $ L_1 + L_2 $ is
dense in $ H $. The latter condition defines a measurable set of pairs
$ (L_1,L_2) $ due to Lemma \ref{6.5}. The former condition may be
expressed in terms of the infinite matrix
\[
M (L_1,L_2) = \( \ip{ f_k(L_1) }{ f_l(L_2) } \)_{k,l}
\]
where $ f_n $ are as in \eqref{6.2}. The relevant set of matrices is
Borel measurable.
\end{proof}

We turn to \sif s. Let $ (\Om,\F,\P) $ be a measure type space and $
\bA (\F) $ the set of all sub-\sif s of $ \F $.\footnote{%
 As before, each \sif\ must contain all $ \P $-negligible sets.}
Then $ \F_0 = \F \bmod 0 $ is a complete Boolean algebra and a Polish
space (recall Footnote 8), and $ \bA (\F) $ may be identified with the
set $ \bA (\F_0) $ of all closed subalgebras of $ \F_0 $. Thus, $ \bA
(\F) = \bA (\F_0) \subset \bF (\F_0) $, and we get measurable maps $
f_n : \bA (\F_0) \to \F_0 $ such that every $ \A \in \bA (\F_0) $ is
generated by (and even the closure of) $ \{ f_1(\A), f_2(\A), \dotsc
\} $.

Striving to a counterpart of \eqref{6.2}, recall the space $ L_0
(\Om,\F,\P) $ of all equivalence classes of measurable maps $ \Om \to
\R $. Every $ X \in L_0 (\Om,\F,\P) $ generates a \sif\ $ \si (X) \in
\bA (\F_0) $. Every $ \A \in \bA (\F_0) $ is $ \si (X) $ for some $ X
\in L_0 (\Om,\F,\P) $; the set of all such $ X $ (for a given $ \A $)
is usually large and non-closed. Nevertheless a selection is
constructed below.

\begin{lemma}
There exists a Borel map $ \A \to X_\A $ from $ \bA (\F_0) $ to $ L_0
(\Om,\F,\P) $ such that $ \si (X_\A) = \A $ for all $ A \in \bA (\F_0)
$.
\end{lemma}

\begin{proof}
One may take
\[
X_\A (\om) = \sum_{n=1}^\infty \frac{ 2 }{ 3^n } \One_{f_n(\A)} (\om)
\, ;
\]
here $ \One_{f_n(\A)} (\om) $ is equal to $ 1 $ if $ \om \in f_n(\A) $
and $ 0 $ otherwise.
\end{proof}

A \sif\ $ \A $ is nonatomic if and only if $ X_\A $ is nonatomic (that
is, $ X_\A^{-1} ( \{ x \} ) $ is negligible for every $ x \in \R
$). Nonatomic elements of $ L_0 (\Om,\F,\P) $ are a Borel
set. Therefore, nonatomic \sif s are a Borel set (in $ \bA (\F_0)
$). Similarly, the number of atoms ($ 0,1,2,\dotsc $ or $ \infty $) is
a Borel function of $ \A $, as well as their (ordered) probabilities.

There is a natural map, $ F \mapsto \si (F) $, from $ \bF (\F_0) $ to
$ \bA (\F_0) $; namely, $ \si (F) $ is the \sif\ generated by $ F
$. The map is measurable, which follows from \eqref{6.1} similarly to
\eqref{6.3}.\footnote{%
 Arbitrary combinations of Boolean operations (union, intersection,
 complement) are used here instead of the linear combinations used in
 \eqref{6.2}.}

Proofs of the following Lemmas \ref{6.9}--\ref{6.11} are left to the
reader, since they are similar to \ref{6.4}--\ref{6.6}
respectively. Especially for \ref{6.11} a hint: $\displaystyle \frac{
Q|_\A }{ P|_\A } = \lim_{n\to\infty} \frac{ Q|_{\A_n} }{ P|_{\A_n} } $
where $ \A_n $ is generated by $ f_1 (\A), \dots, f_n (\A) $.

\begin{lemma}\label{6.9}
Sub-\sif s of $ \F $ are a standard Borel space.
\end{lemma}

\begin{lemma}\label{6.10}
The \sif\ $ \si ( \A_1 \cup \A_2 ) $ generated by $ \A_1, \A_2 $ is a
jointly measurable function of sub-\sif s $ \A_1, \A_2 \subset \F $.
\end{lemma}

The set $ \P $ of equivalent probability measures on $ (\Om,\F) $ is
also a standard Borel space, since it is homeomorphic (in fact,
isometric) to a subset of $ L_1 (P) $ (the choice of $ P \in \P $ does
not matter); the whole $ L_1 (P) $ is a Polish space, and the subset,
consisting of all strictly positive functions whose integral is equal
to $ 1 $, is a Borel subset.

\begin{lemma}\label{6.11}
The Radon-Nikodym density
\[
\frac{ Q|_\A }{ P|_\A }
\]
(treated as an element of $ L_2 (\Om,\F,\P) $ that belongs in fact to
$ L_0 (\Om,\A,\P) $)
is jointly measurable in $ P,Q \in \P $ and $ \A \in \bA (\F) $.
\end{lemma}

Given $ P \in \P $ and $ X \in L_0 (\Om,\F,\P) $, we get a probability
measure on $ \R $, namely, $ S \mapsto P \( \{ \om : X(\om) \in S \}
\) $ for Borel sets $ S \subset \R $. The measure will be called the
distribution of $ X $ w.r.t.\ $ P $ and denoted by $ X(P) $. Note that
$ X(P) $ is jointly measurable in $ X $ and $ P $, that is, $ (X,P)
\mapsto X(P) $ is a Borel map from $ L_0 (\P) \times \P $ to the space
of probability distributions on $ \R $. Indeed, if $ \phi $ is a
bounded continuous function $ \R \to \R $ then $ \int \phi \, d(X(P))
= \int \phi (X(\cdot)) \, dP $ is continuous in $ (X,P) $.

\begin{lemma}\label{6.12}
(a) The set $ \D $ of all pairs $ (\A_1,\A_2) $ such that $ \A_1
\otimes \A_2 $ is well-defined,\footnote{%
 It means existence of $ P \in \P $ that makes $ \A_1, \A_2 $
 independent (recall Def.~\ref{3.3}). When defined, $ \A_1 \otimes
 \A_2 $ is just $ \si (\A_1 \cup \A_2) $.}
is a Borel subset of $ \bA (\F) \times \bA (\F) $, and the map $
(\A_1,\A_2) \mapsto \A_1 \otimes \A_2 $ from $ \D $ to $ \bA (\F) $
is Borel measurable.

(b) For every $ P \in \P $, the map
\[
(\A_1,\A_2) \mapsto \frac{ P|_{\A_1\otimes\A_2} }{ P|_{\A_1} \otimes
P|_{\A_2} } 
\]
from $ \D $ to $ L_0 (\P) $ is Borel measurable.
\end{lemma}

\begin{proof}
(a) The condition $ (\A_1,\A_2) \in \D $ may be expressed in terms of
the measure
\[
\mu_{\A_1,\A_2} = ( X_{\A_1}, X_{\A_2} ) (P)
\]
on $ \R^2 $. (The choice of $ P \in \P $ does not matter.) That is the
joint probability distribution of random
variables $ X_{\A_1}, X_{\A_2} $ on $ (\Om,\F,P) $, and its marginal
distributions are $ X_{\A_1} (P), X_{\A_2} (P) $. Clearly, $
\mu_{\A_1,\A_2} $ is a product measure if and only if $ \A_1, \A_2 $
are independent w.r.t.\ $ P $. Thus, the relevant condition on $
\mu_{\A_1,\A_2} $ says that $ \mu_{\A_1,\A_2} $ must be equivalent to
a product measure. The set of all such measures on $ \R^2 $ is Borel
measurable. The map $ (\A_1,\A_2) \mapsto \A_1 \otimes \A_2 $ is the
restriction to $ \D $ of the map $ (\A_1,\A_2) \mapsto
\si(\A_1\cup\A_2) $ measurable by Lemma \ref{6.10}.

(b) For such $ \A_1, \A_2 $ the map $ (X_{\A_1}, X_{\A_2}) $ is an
isomorphism ($ \bmod 0 $) between $ (\Om,\A_1\otimes\A_2,P) $ and $
(\R^2,\dots,\mu) $, where $ \mu = \mu_{\A_1,\A_2} $ (and the \sif\ of
$ \mu $-measurable subsets of $ \R^2 $ is suppressed in the
notation). Therefore,
\[
\frac{ P }{ P|_{\A_1} \otimes P|_{\A_2} } (\om) = \frac{ \mu }{ \mu_1
\otimes \mu_2 } \( X_{\A_1} (\om), X_{\A_2} (\om) \) \, ,
\]
where $ \mu_1 = \mu_{\A_1} = X_{\A_1} (P) $, $ \mu_2 = \mu_{\A_2} =
X_{\A_2} (P) $ are the marginals of $ \mu $.

I assume in addition that the \sif s $ \A_1, \A_2 $ are nonatomic
(atoms are left to the reader); then an additional transformation $ \R
\to \R $ turns $ \mu_1, \mu_2 $ into Lebesgue measure on $ (0,1)
$.\footnote{%
 One may use the cumulative distribution function $ F_\A (x) = P \( \{
 \om : X_\A (\om) \le x \} \) $; it is continuous (due to the
 nonatomicity), and the random variable $ \om \mapsto F_\A ( X_\A(\om)
 ) $ is distributed uniformly on $ (0,1) $.}
The density of $ \mu $ (w.r.t.\ Lebesgue measure $ \mu_1 \otimes \mu_2
$ on $ (0,1) \times (0,1) $) is a Borel function of $ \mu $ (take an
increasing (refining) sequence of finite partitions of $ (0,1) \times
(0,1) $). The following lemma (or rather, its straightforward
two-dimensional generalization) completes the proof.
\end{proof}

\begin{lemma}
Let $ (\Om,\F,P) $ be a nonatomic probability space. Consider the set
$ \U \subset L_0 (\Om,\F,P) $ of all random variables $ U : \Om \to \R
$ distributed uniformly on $ (0,1) $ (in other words, measure
preserving transformations from $ (\Om,\F,P) $ to $ (0,1) $ with
Lebesgue measure). For each $ U \in \U $ and $ f \in L_0 (0,1) $
consider the composition $ f \circ U $ (that is, $ f(U(\cdot)) $) as
an element of $ L_0 (\Om,\F,P) $. Then

(a) $ U $ is a Borel measurable subset of $ L_0 (\Om,\F,P) $;

(b) $ f \circ U $ is jointly Borel measurable in $ f \in L_0 (0,1) $
and $ U \in \U $.
\end{lemma}

\begin{proof}
(a) Follows immediately from measurability of $ U(P) $ in $ U $.

(b) Let $ Q $ be another probability measure on $ (\Om,\F) $,
equivalent to $ P $. Consider the distribution
\[
( f \circ U ) (Q) = f ( U(Q) ) \, ;
\]
we know that $ f(\mu) $ is jointly measurable in $ f $ and $ \mu $,
and $ U(Q) $ is measurable in $ U $, therefore $ (f \circ U) (Q) $ is
jointly measurable in $ f $ and $ U $. However, distributions $ X(Q) $
for all $ Q $ determine $ X \in L_0 (P) $ uniquely, and moreover, they
generate the Borel \sif\ on $ L_0 (P) $.
\end{proof}

\begin{proposition}\label{6.14}
Let $ (\Om,\F,\P) $ be a measure type space, and $ \cH = \{ (\A,\psi) : \A
\in \bA (\F), \psi \in L_2 (\Om,\A,\P) \} $ the disjoint union of Hilbert
spaces $ L_2 (\Om,\A,\P) $ over all sub-\sif s $ \A \subset \F $. Then

(a) $ \cH $ is (naturally) a standard Borel space;

(b) the set $ \D_1 $ of all pairs $ \( (\A_1,\psi_1), (\A_2,\psi_2) \) \in
\cH \times \cH $ satisfying $ \A_1 = \A_2 $ is a Borel subset of $ \cH
\times \cH $; the map $ \( (\A,\psi_1), (\A,\psi_2) \) \mapsto (\A,\psi_1+\psi_2) $
from $ \D_1 $ to $ \cH $ is Borel measurable; the map $ \( \al, (\A,\psi)
\) \mapsto (\A,\al \psi) $ from $ \R \times \cH $ to $ \cH $ is Borel
measurable; and the map $ \( (\A,\psi_1), (\A,\psi_2) \) \mapsto \ip{ \psi_1 }{
\psi_2 } $ from $ \D_1 $ to $ \R $ is Borel measurable;

(c) the set $ \D_2 $ of all pairs $ \( (\A_1,\psi_1), (\A_2,\psi_2) \) \in
\cH \times \cH $ such that $ \A_1 \otimes \A_2 $ is
well-defined\footnote{%
 See Lemma \ref{6.12}.}
is a Borel subset of $ \cH \times \cH $; and the map $ \(
(\A_1,\psi_1), (\A_2,\psi_2) \) \mapsto \( \A_1 \otimes \A_2, \psi_1 \otimes
\psi_2 \) $ from $ \D_2 $ to $ \cH $ is Borel measurable.
\end{proposition}

\begin{proof}
We choose some $ P \in \P $ and replace each $ L_2 (\Om,\A,\P) $ with
the corresponding $ L_2 (\Om,\A,P) $ according to their unitary
correspondence determined by $ P $,
\[
L_2 (\Om,\A,\P) \ni \psi \mapsto \frac{ \psi }{ \sqrt {P|_\A} } \in
L_2 (\Om,\A,P) \, .
\]
The disjoint union of $ L_2 (\Om,\A,P) $ over all $ \A $ is naturally
a standard Borel space, since it is a Borel subset of the disjoint
union of all subspaces of $ L_2 (\Om,\F,P) $.\footnote{%
 You see, $ L_2 (\Om,\A,\P) $ is not (naturally identified with) a
 subspace of $ L_2 (\Om,\F,\P) $. However, $ L_2 (\Om,\A,P) $ is a
 subspace of $ L_2 (\Om,\F,P) $. Thus, all $ L_2 (\Om,\A,\P) $ become
 embedded into $ L_2 (\Om,\F,\P) $, but the embedding depends on $ P
 $. It may be written as $\displaystyle L_2 (\Om,\A,\P) \ni \psi
 \mapsto \frac{ \psi }{ \sqrt{P|_\A} } \cdot \sqrt P \in L_2 (\Om,\F,P)
 $.}
Thus, a Borel structure
appears also on the disjoint union of $ L_2 (\Om,\A,\P) $, and the
Borel structure does not depend on the choice of $ P $ (since
$\displaystyle \frac{ Q|_\A }{ P|_{\A} } $ is measurable in $ \A $,
see Lemma \ref{6.11}). Item (b) is easily transferred from the
disjoint union of all subspaces of $ L_2 (\Om,\F,P) $ to the disjoint
union of $ L_2 (\Om,\A,\P) $. It remains to prove (c).

Lemma \ref{6.12}(a) gives us measurability of the set $ \D_2 $, and
measurability of $ \A_1 \otimes \A_2 $ in $ (\A_1,\A_2) $. It remains
to verify measurability of $ \psi_1 \otimes \psi_2 $. We know (recall
Sect.~4) that
\[
\psi_1 \otimes \psi_2 = \bigg( \frac{ \psi_1 }{ \sqrt{Q_1} } \otimes
\frac{ \psi_2 }{ \sqrt{Q_2} } \bigg) \cdot \sqrt{ Q_1 \otimes Q_2 }
\in L_2 (\Om,\A_1\otimes\A_2,\P) \, ;
\]
that is, if $ Q \in \P $ is a measure that makes $ \A_1, \A_2 $
independent, then
\[
\frac{ \psi_1 \otimes \psi_2 }{ \sqrt{Q|_{\A_1\otimes\A_2}} } = \frac{
\psi_1 }{ \sqrt{Q|_{\A_1}} } \cdot \frac{ \psi_2 }{ \sqrt{Q|_{\A_2}} }
\quad \text{for } \psi_1 \in L_2 (\Om,\A_1,\P), \, \psi_2 \in L_2
(\Om,\A_2,\P) \, ;
\]
the product in the right-hand side is just a pointwise product of two
functions.\footnote{%
 It is, at the same time, their tensor product, since they are
 independent (w.r.t.\ $ Q $).}
Therefore
\[
\frac{ \psi_1 \otimes \psi_2 }{ \sqrt{ P|_{\A_1\otimes\A_2} } } =
\sqrt{ \frac{ Q|_{\A_1\otimes\A_2} }{ P|_{\A_1\otimes\A_2} } } \cdot
\sqrt{ \frac{ P|_{\A_1} }{ Q|_{\A_1} } } \cdot \sqrt{ \frac{ P|_{\A_2}
}{ Q|_{\A_2} } } \cdot \frac{ \psi_1 }{ \sqrt{ P|_{\A_1} } } \cdot
\frac{ \psi_2 }{ \sqrt{ P|_{\A_2} } } \, .
\]
However, the independence of $ \A_1 $ and $ \A_2 $ under $ Q $ means
that
\[
\frac{ P|_{\A_1\otimes\A_2} }{ Q|_{\A_1\otimes\A_2} } = \frac{
P|_{\A_1\otimes\A_2} }{ Q|_{\A_1} \otimes Q|_{\A_2} } = \frac{
P|_{\A_1\otimes\A_2} }{ P|_{\A_1} \otimes P|_{\A_2} } \cdot \frac{
P|_{\A_1} }{ Q|_{\A_1} } \cdot \frac{ P|_{\A_2} }{ Q|_{\A_2} } \, ,
\]
so,
\[
\frac{ \psi_1 \otimes \psi_2 }{ \sqrt{ P|_{\A_1\otimes\A_2} } } =
\sqrt{ \frac{ P|_{\A_1} \otimes P|_{\A_2} }{ P|_{\A_1\otimes\A_2} } }
\cdot \frac{ \psi_1 }{ \sqrt{P|_{\A_1}} } \cdot \frac{ \psi_2 }{
\sqrt{P|_{\A_2}} } \, .
\]
By Lemma \ref{6.12}(b), $ \dfrac{ P|_{\A_1\otimes\A_2} }{ P|_{\A_1}
\otimes P|_{\A_2} } $ is Borel measurable in $ (\A_1,\A_2) $,
therefore $\displaystyle \frac{ \psi_1 \otimes \psi_2 }{ \sqrt{
P|_{\A_1\otimes\A_2} } } $ is measurable in $\displaystyle \bigg(
\A_1, \A_2, \frac{ \psi_1 }{ \sqrt{P|_{\A_1}} }, \frac{ \psi_2 }{
\sqrt{P|_{\A_2}} } \bigg) $, which means that $ \psi_1 \otimes \psi_2
$ is measurable in $ \( (\A_1,\psi_1), (\A_2,\psi_2) \) $.
\end{proof}

For any Hilbert spaces $ H_1, H_2 $ denote by $ \I (H_1,H_2) $ the set
of all isomorphisms between $ H_1 $ and $ H_2 $, that is, linear
isometric invertible maps $ H_1 \to H_2 $ (of course, $ \I (H_1,H_2) $
is empty if $ H_1, H_2 $ are of different dimension). For any Hilbert
space $ H $ denote by $ \II (H) $ the disjoint union of sets $ \I
(L_1,L_2) $ over all subspaces $ L_1,L_2 \in \bL (H) $. That is, $ \II
(H) $ consists of all triples $ (L_1,L_2,U) $ where $ L_1,L_2 \subset
H $ are subspaces and $ U : L_1 \to L_2 $ is an isomorphism. However,
we may identify each $ U $ with its graph, and $ L_1,L_2 $ with
projections of the graph, $ L_1 (U) = \{ x : (x,y) \in U \} $, $ L_2
(U) = \{ y : (x,y) \in U \} $. Now $ \II (H) $ consists of all
subspaces $ U \in \bL ( H \oplus H ) $ such that $ \| x \| = \| y \| $
whenever $ x \in H $, $ y \in H $, $ (x,y) \in U $. Clearly, $ \II (H)
$ is a Borel subset of $ \bL ( H \oplus H ) $, therefore, a standard
Borel space.

\begin{lemma}\label{6.145}
Let $ H $ be a Hilbert space. Then

(a) the set $ \D $ of all pairs $ \( (L_1,L_2,U), x \) $ such that $
(L_1,L_2,U) \in \II (H) $ and $ x \in L_1 $ is a Borel subset of $ \II
(H) \times H $;

(b) the map $ \( (L_1,L_2,U), x \) \mapsto U(x) $ from $ \D $ to $ H $
is Borel measurable.
\end{lemma}

\begin{proof}
Treating $ U $ as a subspace of $ H \oplus H $ we choose measurable
maps $ f_n : \bL ( H \oplus H ) \to H \oplus H $ such that every $ U $
is spanned by $ \{ f_1(U), f_2(U), \dotsc \} $. We have $ f_n (U) = (
g_n(U), h_h(U) ) $ where $ g_n(U) \in H $, $ h_n(U) \in H $. Applying
the orthogonalization process we ensure that $ g_n (U) $ form an
orthogonal basis of $ L_1 (U) $. Introducing Borel functions $
c_n(x,U) = \ip{ x }{ g_n(U) } $ for $ x \in H $ we have
\begin{gather*}
(U,x) \in \D \equiv \sum_n | c_n(x,U) |^2 = \| x \|^2 \, , \\
(U,x) \in \D \imply U(x) = \sum_n c_n(x,U) h_n(U) \, .
\end{gather*}
\end{proof}

For any FHS space $ G $ we define $ \II (G) $ as consisting
of all triples $ (L_1,L_2,U) $ where $ L_1,L_2 \subset G $ are
subspaces and $ U : L_1 \to L_2 $ is an
FHS-isomorphism. Alternatively, $ \II (G) $ may be thought of as a
subset of $ \bL ( H \oplus H ) $; we'll see (Lemma \ref{6.17}) that it
is a Borel subset, therefore, a standard Borel space.

For any measure type space $ (\Om,\F,\P) $ we define $ \II (\F) $ as
consisting of all triples $ (\A_1,\A_2,U) $ where $ \A_1,\A_2 \in \bA
( \F \bmod 0 ) $ and $ U : \A_1 \to \A_2 $ is an isomorphism of
complete Boolean algebras (it means $ \bmod 0 $ isomorphism between
quotient spaces $ \Om/\A_1 $ and $ \Om/\A_2 $, provided that $
(\Om,\A,\P) $ is a Lebesgue-Rokhlin space). The graph of $ U $ is a
subset of $ ( \F \bmod 0 ) \oplus ( \F \bmod 0 ) = ( \F \oplus \F )
\bmod 0 $, where $ \F \oplus \F $ is the natural \sif\ on the union $
\Om \uplus \Om $ of two disjoint copies of $ \Om $. We identify $ U $
with its graph; it is not just a subset but a \sif; so, $ \II (\F)
\subset \bA ( \F \oplus \F ) $.

\begin{lemma}\label{6.15}
$ \II (\F) $ is a Borel subset of $ \bA ( \F \oplus \F ) $.
\end{lemma}

\begin{proof}
A sub-\sif\ $ \B \subset \F \oplus \F $ belongs to $ \II (\F) $ if and
only if $ \forall \eps \exists \de \forall (A,B) \in \B \( P(A) < \de
\imply P(B) \le \eps ) $; here an element of $ \F \oplus \F $ is
identified with a pair $ (A,B) $ of elements of $ \F $. (The choice of
a measure $ P \in \P $ does not matter.) For such $ \eps $ and $ \de
$, $ \B $ must be disjoint to the \emph{open} set $ \{ (A,B) : P(A) <
\de, P(B) > \eps \} $.
\end{proof}

Let $ (\Om,\F,\P) $ be a measure type space. For any closed set $ F
\subset L_0 (\Om,\F,\P) $ denote by $ \si (F) $ the sub-\sif\
generated by $ F $, that is, the least \sif\ $ \F' \subset \F $ such
that $ F \subset L_0 (\Om,\F',\P) $.

\begin{lemma}\label{6.16}
The map $ F \mapsto \si (F) $ from $ \bF \( L_0 (\Om,\F,\P) \) $ to $
\bA (\Om,\F,\P) $ is Borel measurable.
\end{lemma}

\begin{proof}
Due to \eqref{6.1} it suffices to prove measurability of the map $ f
\mapsto \si (f) $ from $ L_0 (\Om,\F,\P) $ to $ \bA (\Om,\F,\P) $; of
course, $ \si (f) $ means $ \si (\{f\}) $. We know (see the proof of
Lemma \ref{3.8}) that $ f_n \to f $ implies $ \si (f) \subset \liminf
\si (f_n) $. Therefore, for any open set $ V \subset \bA (\Om,\F,\P) $
the set of all $ f $ such that $ \si (f) \cap V = \emptyset $ is
closed.
\end{proof}

\begin{lemma}\label{6.17}
Let $ G $ be an FHS space, then $ \II (G) $ is a Borel subset of $ \bL
( G \oplus G ) $.
\end{lemma}

\begin{proof}
We take a Gaussian type space $ (\Om,\F,\P,G) $\footnote{%
 Suppressing in the notation the distinction between $ G $ and $ G /
 \Const $.}
and consider $ \F \oplus \F $, so that $ G \oplus G \subset L_0 ( \F
\oplus \F ) $. Every subspace $ U \in \bL ( G \oplus G ) $ generates
a sub-\sif\ $ \si(U) \in \bA ( \F \oplus \F ) $. It is easy to see
that $ U \in \II (G) $ if and only if $ \si(U) \in \II(\F) $. However,
$ \si(U) $ is measurable in $ U $ by Lemma \ref{6.16}, and $ \II (\F)
$ is measurable by Lemma \ref{6.15}.
\end{proof}

For any measure type space $ (\Om,\F,\P) $ we define $ \II (\P)
$\footnote{%
 Sorry for the clumsy notation: $ \II(\F) $ for \sif s, but $ \II(\P)
 $ for square roots of measures.}
as consisting of all triples $ (\A_1,\A_2,U) $ where $ \A_1,\A_2
\subset \F $ are sub-\sif s and $ U \in \I \( L_2(\Om,\A_1,\P),
\linebreak[0] L_2(\Om,\A_2,\P) \) $ is a linear isometry. Spaces $ L_2
(\Om,\A,\P) $ are not naturally embedded into $ L_2 (\Om,\F,\P) $;
however, we may choose some measure $ P \in \P $ and embed all $ L_2
(\Om,\A,\P) $ into $ L_2 (\Om,\F,P) $ by
\[
L_2 (\Om,\A,\P) \ni \psi \mapsto \frac{ \psi }{ \sqrt{P|_\A} } \in L_2
(\Om,\A,P) \subset L_2 (\Om,\F,P) \, .
\]
We have a bijective correspondence between $ \II (\P) $ and $ \II \(
L_2 (\Om,\F,P) \) $, which turns $ \II (\P) $ into a standard Borel
space. Its Borel structure does not depend on the choice of $ P \in \P
$ (recall Lemma \ref{6.11}), but the correspondence depends on $ P $.

Take an element $ (\A_1,\A_2,U) $ of $ \II (\F) $ (this time we prefer
a triple to a graph). The isomorphism $ U $ between \sif s $ \A_1,
\A_2 $ induces naturally an isomorphism (linear isometry) between
Hilbert spaces $ L_2 (\Om,\A_1,\P) $ and $ L_2 (\Om,\A_2,\P)
$. Namely, if measures $ P_1 \in \P|_{\A_1} $, $ P_2 \in \P|_{\A_2} $
satisfy $ P_2 ( U(A) ) = P_1 (A) $ for all $ A \in \A_1 $, then the
vector $ \psi_2 = \sqrt{P_2} \in L_2 (\Om,\A_2,\P) $ corresponds to
the vector $ \psi_1 = \sqrt{P_1} \in L_2 (\Om,\A_1,\P) $. (Such
vectors are not a linear set, but span the Hilbert spaces.) So, we
have a map from $ \II (\F) $ to $ \II (\P) $.

\begin{lemma}\label{6.18}
The map from $ \II (\F) $ to $ \II (\P) $ is Borel measurable.
\end{lemma}

\begin{proof}
Consider some $ U \in \II (\F) $ treated as a sub-\sif\ of $ \F \oplus
\F $. Given some $ P_1 \in \P $, we introduce on $ U $ a measure $ P
(A,B) = P_1 (A) $ for $ (A,B) \in U $; here, as before, an element of
$ \F \oplus \F $ is represented by a pair $ (A,B) $ where $ A,B \in \F
$. There is also a measure $ P_2 $ on the \sif\ $ \A_2 = \{ B : (A,B)
\in U \} $ such that $ P_2 (B) = P (A,B) = P_1 (A) $ whenever $ (A,B)
\in U $. The pair $ \( \sqrt{P_1}, \sqrt{P_2} \) $ belongs to the
graph $ U_1 \in \II (\P) $ that corresponds to $ U $. Such pairs for
all $ P_1 \in \P $ (or for a countable dense subset) span the graph $
U_1 $. It remains to prove that, for a given $ P_1 $ and arbitrary $ U
$, the pair $ \( \sqrt{P_1}, \sqrt{P_2} \) $ is measurable in $ U $
(you see, $ P_2 $ depends implicitly on $ U $). According to our
definition of the Borel structure on $ \II (\P) $, we have to prove
measurability in $ U $ of the density $\displaystyle \frac{ P_2 }{ P_1
|_{\A_2} } $. The latter is the restriction (to the second copy of $
\Om $) of a density on the doubled space, $ \Om \uplus \Om $, namely,
$ \displaystyle \frac{P'|_U}{P''|_U} $, where measures $ P', P'' $ on
$ \F \oplus \F $ are defined (irrespective of $ U $) by $ P'(A,B) =
P_1(A) $, $ P''(A,B) = P_1(B) $. We apply Lemma \ref{6.11} to $ P',
P'' $ on $ \Om \uplus \Om $. Though, $ P' $ and $ P'' $ are not
equivalent, but one can consider, say, $ \displaystyle \frac{ 2P'+P''
}{ P'+P'' } $.
\end{proof}

\begin{proposition}\label{6.19}
Let $ (\Om,\F,\P,G) $ be a Gaussian type space, and $ G_0 = G / \Const
$ the corresponding FHS space. Then the natural map from $ \II (G_0) $
to $ \II (\P) $ is Borel measurable.
\end{proposition}

\begin{proof}
We know that every $ U \in \II (G) $ generates $ \si(U) \in \II(\F) $
(see the proof of Lemma \ref{6.17}), and the map $ U \mapsto \si(U) $
is measurable. The transition from $ \si(U) $ to the corresponding
element of $ \II(\P) $ is measurable by Lemma \ref{6.18}.
\end{proof}

%% file: sect7.tex
\begin{definition}\label{7.1}
A \emph{sum system} consists of a two-parameter family $ (G_{a,b}) $
of FHS-spaces $ G_{a,b} $, given for $ -\infty < a < b < +\infty $,
embedded into a single linear space $ G_{-\infty,+\infty} $, and a
one-parameter group $ (U_t) $ of linear maps $ U_t :
G_{-\infty,+\infty} \to G_{-\infty,+\infty} $ for $ t \in \R $, such
that

(a) $ G_{a,c} = G_{a,b} \oplus G_{b,c} $ (in the FHS sense) whenever $
-\infty < a < b < c < +\infty $;\footnote{%
 Thus, $ H_{a,b} $ and $ H_{b,c} $ must be linear subspaces of $
 H_{a,c} $; and their FHS structures must be inherited from $ H_{a,c}
 $; and they must be orthogonal in some admissible norm. Note that the
 norm may depend on $ b $.}

(b) $ U_t : G_{a,b} \to G_{a+t,b+t} $ is an isomorphism of FHS spaces, 
whenever $ -\infty < a < b < +\infty $ and $ -\infty < t < +\infty
$.\footnote{%
 Thus, $ U_t $ must map $ G_{a,b} $ onto $ G_{a+t,b+t} $ and send an
 admissible norm into an admissible norm.}

(c) $ (U_t) $ is strongly continuous in the sense that $ \| U_t x - x
\| \to 0 $ when $ t \to 0 $, whenever $ a < b < c < d $, $ x \in
G_{b,c} $, and the norm is taken in $ G_{a,b} $ (which is correct for
$ t $ small enough).
\end{definition}

The structure is local in the sense that the global space $
G_{-\infty,+\infty} $ is equipped with a linear structure only, not an
FHS structure, nor even a topology. One may assume that $
G_{-\infty,+\infty} $ is just the union of all $ G_{a,b} $, since
anyway, only the local spaces $ G_{a,b} $ will be used.

Given a sum system $ \( (G_{a,b}), (U_t) \) $, we may introduce (as
explained in Sect.\ 5) Hilbert spaces $ H_{a,b} = \Exp (G_{a,b}) $
satisfying (under the usual identification)
\[
H_{a,b} \otimes H_{b,c} = H_{a,c}
\]
whenever $ -\infty < a < b < c < +\infty $. Given $ a < b < c < d $, $
\psi_1 \in H_{a,b} $, $ \psi_2 \in H_{b,c} $, $ \psi_3 \in H_{c,d} $,
we may calculate $ \psi_1 \otimes \psi_2 \otimes \psi_3 $ as $ (
\psi_1 \otimes \psi_2 ) \otimes \psi_3 $ or $ \psi_1 \otimes ( \psi_2
\otimes \psi_3 ) $, which is the same (recall Lemma \ref{3.5}). The
two-parameter families may be reduced to one-parameter families using
$ (U_t) $; namely, for $ s,t \in (0,\infty) $,
\begin{equation}
\begin{aligned}
G_{0,s+t} &= G_{0,s} \oplus U_s G_{0,t} \, , \\
H_{0,s+t} &= H_{0,s} \otimes \( \Exp ( U_s|_{G_{0,t}} ) \) H_{0,t} \,
,
\end{aligned}
\end{equation}
where $ \( \Exp ( U_s|_{G_{0,t}} ) \) : H_{0,t} \to H_{s,s+t} $ is the
unitary operator corresponding to the FHS isomorphism $ U_s
|_{H_{0,t}} : H_{0,t} \to H_{s,s+t} $. The binary operation of tensor
product, $ (s,\psi_1), (t,\psi_2) \mapsto \( s+t, \psi_1 \otimes (
\Exp (U_s|_{G_{0,t}}) ) \psi_2 \) $ (for $ \psi_1 \in H_{0,s} $, $
\psi_2 \in H_{0,t} $) is associative. That is the algebraic part of
the `product system' structure defined by W.~Arveson
\cite[Def.~1.4]{Ar89}. It is not the whole story, since some
measurability in $ s,t $ is needed. Namely, the disjoint union of
spaces $ H_{0,s} $ must be a standard Borel space, and the tensor
product must be a measurable binary operation.

The disjoint union is the set of all pairs $ (s,\psi) $ such that $ s
\in (0,\infty) $ and $ \psi \in H_{0,s} $. We take some $ T \in
(0,\infty) $; it is enough to consider $ s \in (0,T) $ rather than $
(0,\infty) $. We have a Gaussian type space $ (\Om,\F,\P,G) $ such
that $ G / \Const = G_{0,T} $, and sub-\sif s $ \F_{a,b} \subset \F $
such that $ L_2 (\Om,\F_{a,b},\P) = H_{a,b} $ whenever $ (a,b) \subset
(0,T) $.

\begin{lemma}\label{7.3}
The map $ (a,b) \mapsto G_{a,b} $ from the triangle $ \{ (a,b) : 0 \le
a < b \le T \} $ to the Borel space $ \bL (G_{0,T}) $ of subspaces of
$ G_{0,T} $ is Borel measurable.
\end{lemma}

\begin{proof}
It is enough to consider the case $ 0 \le a < t < b \le T $ for an
arbitrary $ t \in (0,T) $. The equality $ G_{a,b} = G_{a,t} \oplus
G_{t,b} $, in combination with Lemma \ref{6.5}, reduces the problem to
measurability of $ G_{a,t} $ in $ a $, and $ G_{t,b} $ in $ b
$. However, a \emph{monotone} function is always measurable.
\end{proof}

So, $ G_{a,b} $ is measurable in $ (a,b) $. It follows by Lemma
\ref{6.16} that $ \F_{a,b} $ is measurable in $ (a,b) $. Proposition
\ref{6.14} gives us a Borel structure on the disjoint union of spaces
$ H_{a,b} $ (over all $ a,b $ satisfying $ 0 \le a < b \le T $); it is
compatible with linear operations and scalar product; and the map
\[
\( ((a,b),\psi_1), ((b,c),\psi_2) \) \mapsto ( (a,c), \psi_1 \otimes
\psi_2 )
\]
is Borel measurable (here $ \psi_1 \in H_{a,b} $, $ \psi_2 \in H_{b,c}
$, $ 0 \le a < b < c \le T $).

Given $ a,b $ such that $ 0 \le a < b \le T $, we may treat the
restriction $ U_{b-a} |_{G_{0,a}} $ as an FHS isomorphism $ G_{0,a}
\to G_{b-a,b} $, therefore an element of the Borel space $ \II
(G_{0,T}) $ introduced in Sect.~6.

\begin{lemma}
The map $ (a,b) \mapsto U_{b-a} |_{G_{0,a}} $ from the triangle $ \{
(a,b) : 0 \le a < b \le T \} $ to $ \II (G_{0,T}) $ is Borel
measurable.
\end{lemma}

\begin{proof}
By Lemma \ref{7.3}, $ G_{0,a} $ is measurable in $ a $. According to
\eqref{6.1} there are $ f_n (a) \in G_{0,a} $, measurable in $ a $,
such that every $ G_{0,a} $ is spanned by $ \{ f_1(a), f_2(a), \dotsc
\} $. Therefore the graph of $ U_{b-a} |_{G_{0,a}} $ is spanned by
pairs $ \( f_n(a), U_{b-a} f_n(a) \) $. Each pair is measurable in $ a
$ and continuous in $ b $ (recall \ref{7.1}(c)), therefore, measurable
in $ (a,b) $ (see \cite[Th.~3.1.30]{Sr}).
\end{proof}

\begin{sloppypar}
So, $ U_{b-a} |_{G_{0,a}} $ is measurable in $ (a,b) $. It follows by
Proposition \ref{6.19} that $ \Exp ( U_{b-a} |_{G_{0,a}} ) $ is also
measurable in $ (a,b) $. Lemma \ref{6.15} shows that $ \Exp ( U_{b-a}
|_{G_{0,a}} ) \psi $ is jointly measurable in $ \psi \in H_{0,a} $ and
$ a,b $. It follows that $ \psi_1 \otimes \( \Exp ( U_{b-a}
|_{G_{0,a}} ) \) \psi_2 $ is jointly measurable in $ a $, $ b $, $
\psi_1 \in H_{b-a} $, $ \psi_2 \in H_a $, which proves the following
result.
\end{sloppypar}

\begin{theorem}\label{7.5}
If $ \( (G_{a,b}), (U_t) \) $ is a sum system, then Hilbert spaces $
H_{a,b} = \Exp (G_{a,b}) $ with the natural identification $ H_{0,s+t}
= H_{0,s} \otimes \( \Exp ( U_s |_{G_{0,t}} ) \) H_{0,t} $ form a
product system.
\end{theorem}

The product system may be called the exponential of the given sum
system.

%% file: sect8.tex
An isomorphism of two product systems is defined \cite[p.~6]{Ar89} as
a family $ (V_t) $ of linear isomorphisms $ V_t : H'_{0,t} \to
H''_{0,t} $ between corresponding Hilbert spaces that respects the two
structures on the disjoint union of the Hilbert spaces, namely, the
binary operation of tensor multiplication, and the Borel
\sif. Assuming that the two product systems are exponentials of two
given sum systems $ \( (G'_{a,b}), (U'_t) \) $ and $ \( (G''_{a,b}),
(U''_t) \) $, we may redefine equivalently an isomorphism as a
two-parameter family $ (V_{a,b}) $ of unitary operators that satisfy
\begin{gather*}
V_{a,b} : H'_{a,b} \to H''_{a,b} \text{ unitarily,} \\
V_{a,c} = V_{a,b} \otimes V_{b,c} \, , \\
\Exp ( U''_t |_{G''_{a,b}} ) V_{a,b} = V_{a+t,b+t} \Exp ( U_t
  |_{G'_{a,b}} ) \, , \\
\end{gather*}
whenever $ -\infty < a < b < c < +\infty $ and $ -\infty < t < +\infty
$ (as before, $ H'_{a,b} = \Exp ( G'_{a,b} ) $, $ H''_{a,b} = \Exp (
G''_{a,b} ) $), and respects the Borel structure on the disjoint union
of Hilbert spaces.

For now $ H_E $, as well as $ G_E $ and $ \F_E $, are defined only
when $ E $ is an interval, $ E = (a,b) $. However, they may be defined
for any elementary set $ E $, that is, a union of a finite number of
intervals. Given $ -\infty < a < b < c < d < +\infty $, we have
\[
\underbrace{ H_{a,d} }_{ H_{(a,d)} } = H_{a,b} \otimes H_{b,c} \otimes
H_{c,d} = \underbrace{ ( H_{a,b} \otimes H_{c,d} ) }_{
H_{(a,b)\cup(c,d)} } \otimes \underbrace{ H_{b,c} }_{ H_{(b,c)} } \, .
\]
The same for any finite number of intervals. We get $ H_{E_1\cup E_2}
= H_{E_1} \otimes H_{E_2} $ when $ E_1 \cap E_2 = \emptyset $. Dealing
with elementary sets we neglect boundary points, treating, say, $
(a,b) \cup (b,c) $ as $ (a,c) $. Also, $ H_E = H_{E_1} \otimes \dots
\otimes H_{E_n} $ whenever $ E_1, \dots, E_n $ are pairwise disjoint
and $ E = E_1 \cup \dots \cup E_n $. Similarly, $ G_E = G_{E_1} \oplus
\dots \oplus G_{E_n} $ (in the FHS sense), and $ \F_E = \F_{E_1}
\otimes \dots \otimes \F_{E_n} $ (recall \ref{3.2}--\ref{3.5}).

\begin{proposition}\label{8.1}
Let two sum systems $ \( (G'_{a,b}), (U'_t) \) $ and $ \( (G''_{a,b}),
(U''_t) \) $ be such that the corresponding product systems are
isomorphic. Let $ E_1, E_2, \dotsc \subset (0,1) $ be elementary
sets. Then the following two conditions are equivalent.

(a) $ \liminf G'_{(0,1)\setminus E_n} = G'_{(0,1)} $ and $ \limsup
G'_{E_n} = \{0\} $;

(b) $ \liminf G''_{(0,1)\setminus E_n} = G''_{(0,1)} $ and $ \limsup
G''_{E_n} = \{0\} $.

\noindent (The FHS spaces are treated as subspaces of $ G'_{(0,1)} $,
$ G''_{(0,1)} $ respectively.)
\end{proposition}

\begin{proof}
Theorem \ref{5.3} allows us to reformulate the conditions in terms of
density matrices in product systems, thus making explicit their
invariance under isomorphisms.
\end{proof}

%% file: sect9.tex
Consider a scalar product of the form
\begin{equation}\label{9.1}
\ip f g = \iint f(s) g(t) B(s-t) \, ds dt
\end{equation}
assuming that $ B : \R \setminus \{0\} \to \R $ is continuous outside
of the origin, and $ B(-t) = B(t) $ for all $ t \in (0,\infty) $, and
$ \int_0^1 |B(t)| \, dt < \infty $. The scalar product is well-defined
whenever $ f,g \in L_2(-M,M) $, $ M \in (0,\infty) $. We assume that $
B $ is positively definite in the sense that $ \ip f f \ge 0 $ for all
such $ f $.

Denote by $ G_{a,b} $ the completion of $ L_2 (a,b) $ w.r.t.\ the
scalar product \eqref{9.1}, then $ G_{a,b} $ is a Hilbert space. We
introduce operators $ U_t $ by $ (U_t f) (s) = f (s-t) $ for $ f \in
L_2 (a,b) $ and extend $
U_t $ by continuity to any $ G_{a,b} $; thus, $ U_t : G_{a,b} \to
G_{a+t,b+t} $ is a unitary operator, and $ U_s U_t = U_{s+t} $, and $
\| U_t f - f \|_{G_{a,d}} \to 0 $ for $ t \to 0 $, if $ f \in G_{b,c}
$ and $ a < b < c < d $ (since it holds for the dense subset of
continuous functions $ f $).

In order to get a sum system (as defined by \ref{7.1}) we need to
ensure that $ G_{a,c} = G_{a,b} \oplus G_{b,c} $ in the FHS sense. The
property will be verified for $ B $ such that
\begin{equation}\label{9.2}
\begin{gathered}
\exists \eps > 0 \; \forall t \in (0,\eps) \; B(t) = \frac1{ t \ln^\al
  (1/t) } \, ; \\
\text{on $ (0,\infty) $ the function $ B(\cdot) $ is positive,
  decreasing and convex;}
\end{gathered}
\end{equation}
here $ \al \in (1,\infty) $ is a parameter. Such $ B $ is positively
definite, since it is an integral combination (with positive weights)
of `triangle' functions of the form $ t \mapsto \max ( 0, a-|t| ) $,
and maybe a positive constant function.

We consider $ G_{-T,0} $ and $ G_{0,T} $ for an arbitrary $ T \in
(0,\infty) $. To this end we introduce $ X_k \in G_{0,T} $ and $ Y_k
\in G_{-T,0} $ by
\begin{equation}\label{9.3}
\begin{split}
X_k (t) &= \One_{(0,T)} (t) \cdot \exp ( 2\pi i k t / T ) \, , \\
Y_k (t) &= \One_{(-T,0)} (t) \cdot \exp ( 2\pi i k t / T )
\end{split}
\end{equation}
for $ k \in \Z $; of course, $ \One_{(a,b)} $ is the indicator of $
(a,b) $. Clearly, $ X_k $ span $ G_{0,T} $, and $ Y_k $ span $
G_{-T,0} $. An elementary calculation gives
\begin{equation}\label{9.4}
\begin{gathered}
\ip{ X_k }{ X_k } = \ip{ Y_k }{ Y_k } = 2 \int_0^T (T-t) B(t) \cos
  (2\pi k t / T ) \, dt \, , \\
\ip{ X_k }{ X_l } = \ip{ Y_k }{ Y_l } = - \frac{ T }{ \pi (k-l) }
  \int_0^T B(t) \( \sin ( 2\pi k t / T ) - \sin ( 2\pi l t / T ) \) \,
  dt \, , \\
\ip{ X_k }{ Y_k } = \int_0^{2T} \min (t,2T-t) B(t) \exp ( 2\pi i k t /
  T ) \, dt \, , \\
\ip{ X_k }{ Y_l } = - i \frac{ T }{ 2\pi (k-l) } \int_0^{2T} B(t)
  \sgn(T-t) \( \exp ( 2\pi i k t / T ) - \exp ( 2\pi i l t / T ) \) \,
  dt
\end{gathered}
\end{equation}
for $ k \ne l $. We want to estimate $ \ip{ X_k }{ X_l } $ and $ \ip{
X_k }{ Y_l } $ from above. These are increments of Fourier transforms,
thus we want to differentiate these Fourier transforms. The
singularity of $ B $ at the origin contributes a term that decays
slowly (near $ \infty $) and is monotone. Jumps outside the origin (at
$ T $ and $ 2T $) contribute terms that decay much faster, but
oscillate. After differentiation, these oscillating terms dominate the
monotone term. However, we need Fourier transforms only on the lattice
$ (2\pi / T) \Z $, thus we have a freedom to change the given
functions without changing their Fourier transforms on the
lattice. We'll use the freedom for eliminating the jumps.

Note that $ \int e^{i\la t} (U_s f) (t) \, dt = e^{i\la s} \int
e^{i\la t} f(t) \, dt $, therefore $ \int e^{i\la t} ( U_T f - f ) (t)
\, dt = 0 $ for $ \la \in (2\pi / T) \Z $. We use a piecewise linear $
f $ for correcting $ B $; namely, we define
\[\begin{split}
b_1 (t) &= \begin{cases}
  B(t) - B(T) \cdot \frac{T-t}{T} & \text{for $ 0 < t \le T $}, \\
  B(T) \cdot \frac{2T-t}{T} & \text{for $ T \le t \le 2T $}, \\
  0  & \text{for other $ t $},
 \end{cases} \\
\hat b_1(\la) &= \int_0^\infty e^{i\la t} b_1(t) \, dt \, ,
\end{split}\]
then $ b_1 $ is continuous on $ (0,\infty) $, and
\[
\ip{ X_k }{ X_l } = \ip{ Y_k }{ Y_l } = -\frac{T}{\pi(k-l)} \Im \(
\hat b_1 (2\pi k/T) - \hat b_1 (2\pi l/T) \)
\]
for $ k \ne l $. Similarly,
\[\begin{split}
b_2 (t) &= \begin{cases}
  B(t) - B(t+T) - (B(T)-B(2T)) \cdot \frac{T-t}{T} & \text{for $ 0 < t
    \le T $}, \\
  ( B(T) - B(2T) ) \cdot \frac{2T-t}{T} & \text{for $ T \le t \le 2T
    $}, \\
  0  & \text{for other $ t $},
 \end{cases} \\
\hat b_2(\la) &= \int_0^\infty e^{i\la t} b_2(t) \, dt \, ; \\
\ip{ X_k }{ Y_l } &= -i \frac{T}{2\pi(k-l)} \( \hat b_2 (2\pi k/T) -
  \hat b_2 (2\pi l/T) \) \quad \text{for } k \ne l \, .
\end{split}\]
It is easy to check that both $ b_1 $ and $ b_2 $ satisfies the
conditions of the following lemma, provided that $ B $ satisfies
\eqref{9.2}.\footnote{%
 Finite variation of $ (tB(t))' $ on any $ [\eps,1/\eps] $ follows
from increase of $ (tB(t))' - 2B(t) $.}

\begin{lemma}\label{9.5}
Assume that $ \al \in (1,\infty) $, and a function $ b : (0,\infty)
\to \R $ is continuous, and the difference $\displaystyle b(t) -
\frac{ \One_{(0,\eps)} (t) }{ t \ln^\al (1/t) } $ is of finite
variation on $ (0,\infty) $ for some (therefore, every) $ \eps \in
(0,1) $, and $ b(t) = 0 $ for all $ t $ large enough. Assume also that
the function $ t \mapsto t b(t) $ is absolutely continuous on $
(0,\infty) $, and the difference $\displaystyle \( t b(t) \)' -
\al \frac{ \One_{(0,\eps)} (t) }{ t \ln^{\al+1} (1/t) } $ is of finite
variation on $ (0,\infty) $ for some (therefore, every) $ \eps \in
(0,1) $. Then the function $ \hat b(\la) = \int_0^\infty e^{i\la t}
b(t) \, dt $ satisfies
\[\begin{split}
\hat b(\la) &= \frac1{\al-1} \frac1{ \ln^{\al-1} \la } + O \Big(
  \frac1{ \ln^\al \la } \Big) \qquad \text{for } \la \to +\infty \, ,
  \\
\frac{d}{d\la} \hat b(\la) &= - \frac1{\la \ln^\al \la} + O \Big(
  \frac1{ \la \ln^{\al+1} \la } \Big) \qquad \text{for } \la \to
  +\infty \, .
\end{split}\]
\end{lemma}

\begin{proof}
Choosing any $ \eps \in (0,1) $ we have for large $ \la $
\begin{multline*}
\hat b(\la) = \underbrace{ \int_0^\infty e^{i\la t} \Big( b(t) -
  \frac{ \One_{(0,\eps)} (t) }{ t \ln^\al (1/t) } \Big) \, dt }_{
  O(1/\la) } + \int_0^\eps e^{i\la t} \frac1{ t \ln^\al (1/t) } \, dt
  = \\
= \bigg( \int_0^{1/\la} + \int_{1/\la}^\eps \bigg) e^{i\la t} \frac1{
  t \ln^\al (1/t) } \, dt + O(1/\la) = \\
= \int_0^{1/\la} \frac1{ t \ln^\al (1/t) } \, dt + \int_0^{1/\la}
  \frac{ e^{i\la t} - 1 }{ t \ln^\al (1/t) } \, dt + \int_{1/\la}^\eps
  e^{i\la t} \frac1{ t \ln^\al (1/t) } \, dt + O(1/\la) \, ;
\end{multline*}
\[\begin{split}
& \int_0^{1/\la} \frac1{ \ln^\al (1/t) } \frac{dt}t = \frac1{\al-1}
  \frac1{ \ln^{\al-1} \la } \, ; \\
& \bigg| \int_0^{1/\la} \frac{ e^{i\la t} - 1 }{ t \ln^\al (1/t) } \,
  dt \bigg| \le \int_0^{1/\la} \frac{ \la t }{ t \ln^\al (1/t) } \, dt
  = \la \int_0^{1/\la} \frac{ dt }{ \ln^\al (1/t) } \le \frac1{
  \ln^\al \la } \, ; \\
& \int_{1/\la}^\eps \frac{ e^{i\la t} }{t} \frac1{\ln^\al (1/t)} \, dt =
  \int_{1/\la}^\eps \frac1{ \ln^\al (1/t)} \, d \( \Ci(\la t) + i
  \Si(\la t) \) \, ,
\end{split}\]
where $\displaystyle \Ci(t) = - \int_t^\infty \frac{ \cos u }{ u } \,
du $, $\displaystyle \Si(t) = - \int_t^\infty \frac{ \sin u }{ u } \,
du $. Taking into account that $ \Ci(t) = O(1/t) $ and $ \Si(t) =
O(1/t) $, we get
\begin{multline*}
\int_{1/\la}^\eps e^{i\la t} \frac1{ t \ln^\al (1/t) } \, dt =
\underbrace{ \frac{ \Ci(\la t) + i\Si(\la t) }{ \ln^\al (1/t) }
  \bigg|_{1/\la}^\eps }_{ O (1/\ln^\al \la) } - \\
- \underbrace{ \int_{1/\la}^\eps \( \Ci(\la t) + i\Si(\la t) \) \cdot
  \frac{ \al }{ t \ln^{\al+1} (1/t) } \, dt }_{ O ( \int_{1/\la}^\eps
  \frac1{\la t} \cdot \frac{ dt }{ t \ln^{\al+1} (1/t) } ) } \, ;
\end{multline*}
\begin{multline*}
\frac1\la \int_{1/\la}^\eps \frac{ dt }{ t^2 \ln^{\al+1} (1/t) } =
  \frac1\la \bigg( \int_{1/\la}^{1/\sqrt\la} + \int_{1/\sqrt\la}^\eps
  \bigg) \frac{ dt }{ t^2 \ln^{\al+1} (1/t) } \le \\
\le \frac1\la \frac1{ \ln^{\al+1} \sqrt\la } \underbrace{
  \int_{1/\la}^\infty \frac{ dt }{ t^2 } }_{ \la } + \frac1\la
  \frac1{\ln^{\al+1} (1/\eps)} \underbrace{ \int_{1/\sqrt\la}^\infty
  \frac{ dt }{ t^2 } }_{ \sqrt\la } = O \Big( \frac1{ \ln^{\al+1} \la
  } \Big) \, .
\end{multline*}
So,
\[
\hat b(\la) = \frac1{\al-1} \frac1{ \ln^{\al-1} \la } + O \Big(
\frac1{ \ln^\al \la } \Big) \quad \text{for } \la \to +\infty \, ,
\]
which is the first claim of the lemma. In order to prove the
second claim we note that the only properties of the function $ b(t) $
used till now are the finite variation of $ b(t) - \frac{
\One_{(0,\eps)} (t) }{ t \ln^\al (1/t) } $, and $ b(t)=0 $ for large $
t $. Therefore the same argument may be applied to the function $
\frac1\al \( t b(t) \)' $ w.r.t.\ $ \al+1 $:
\[
\int_0^\infty e^{i\la t} \frac1\al \( t b(t) \)' \, dt = \frac1\al
\frac1{ \ln^\al \la } + O \Big( \frac1{ \ln^{\al+1} \la } \Big) \, .
\]
Hence
\begin{multline*}
\frac{d}{d\la} \hat b(\la) = \int_0^\infty e^{i\la t} i t b(t) \, dt =
  \frac1\la \int_0^\infty t b(t) ( e^{i\la t} )' \, dt = \\
= - \frac1\la \int_0^\infty \( t b(t) \)' e^{i\la t} \, dt = - \frac1{
  \la \ln^\al \la } + O \Big( \frac1{ \la \ln^{\al+1} \la } \Big) \, .
\end{multline*}
\end{proof}

\begin{lemma}\label{9.6}
Let $ B $ satisfy \eqref{9.2} and $ X_k, Y_k $ be defined by
\eqref{9.3}; then
\[
\sum_{m,n: m \ne n} \frac{ | \ip{ X_m }{ X_n } |^2 }{ \| X_m \|^2 \|
  X_n \|^2 } < \infty \, , \qquad
\sum_{m,n} \frac{ | \ip{ X_m }{ Y_n } |^2 }{ \| X_m \|^2 \| Y_n \|^2 }
  < \infty \, .
\]
\end{lemma}

\begin{proof}
First, the function $ t B(t) $ on $ [0,T] $ is of finite variation,
thus, using \eqref{9.4} and Lemma \ref{9.5},
\begin{multline*}
\ip{ X_k }{ X_k } = 2T \underbrace{ \int_0^T B(t) \cos ( 2\pi k t / T
  ) \, dt }_{\hat b_1 (2\pi k/T)} - 2 \underbrace{ \int_0^T t B(t)
  \cos ( 2\pi k t / T ) \, dt }_{O(1/|k|)} \sim \\
\sim 2T \cdot \frac1{\al-1} \cdot \frac1{ \ln^{\al-1} (2\pi|k|/T) }
  \sim \frac{2T}{\al-1} \frac1{ \ln^{\al-1} |k| } \, ;
\end{multline*}
\[
\frac1{\|X_k\|^2} = O \( \ln^{\al-1} |k| \)
\]
for $ k \to \pm\infty $.

\begin{sloppypar}
Second, $ \ip{ X_k }{ Y_k } = \int_0^{2T} \min(t,2T-t) B(t) \exp(2\pi
i k t/T) = O ( 1/|k| ) $, since $ \min(t,2T-t) B(t) $ is of finite
variation on $ [0,2T] $. Hence $ | \ip{ X_n }{ Y_n } |^2 / ( \| X_n
\|^2 \| Y_n \|^2 ) = O \( \frac1{n^2} \ln^{2\al-2} |n| \) $ and
\[
\sum_n \frac{ | \ip{ X_n }{ Y_n } |^2 }{ \| X_n \|^2 \| Y_n \|^2 } <
\infty \, .
\]
\end{sloppypar}

Third, $ | \ip{ X_n }{ X_{-n} } | = O (1/|n|) $ and $ | \ip{ X_n }{
Y_{-n} } | = O (1/|n|) $, hence
\[
\sum_n \frac{ | \ip{ X_n }{ X_{-n} } |^2 }{ \| X_n \|^2 \| X_{-n} \|^2
} < \infty \, , \qquad
\sum_n \frac{ | \ip{ X_n }{ Y_{-n} } |^2 }{ \| X_n \|^2 \| Y_{-n} \|^2
} < \infty \, .
\]

It is enough to prove that\footnote{%
 You see, $ \ln |n| $ is replaced with $ \ln (|n|+2) $ in order to
 cover the small values, $ n = -1,0,1 $.}
\[
\sum_{m,n:\,m\pm n\ne0} \frac{ \ln^{\al-1} (|m|+2) \cdot \ln^{\al-1}
(|n|+2) }{ (m-n)^2 } | \hat b (2\pi m/T) - \hat b (2\pi n/T) |^2 <
\infty
\]
for every function $ b $ as in Lemma \ref{9.5}. Taking into account
that $ \hat b(-\la) = \overline{ \hat b(\la) } $ we transform it into
\begin{multline*}
\sum_{m,n:\,0\le m<n} \ln^{\al-1} (m+2) \ln^{\al-1} (n+2) \cdot \bigg(
  \frac{ | \hat b(2\pi m/T) - \hat b(2\pi n/T) |^2 }{ (n-m)^2 } + \\
+ \frac{ | \overline{ \hat b(2\pi m/T) } - \hat b(2\pi n/T) |^2 }{
  (n+m)^2 } \bigg) < \infty \, ;
\end{multline*}
\begin{multline*}
\sum_{m,n:\,0\le m<n} \ln^{\al-1} (m+2) \ln^{\al-1} (n+2) \cdot \\
\cdot \bigg( \( \Re \hat b(2\pi m/T) - \Re \hat b(2\pi n/T) \)^2 \cdot
  \Big( \frac1{ (n-m)^2 } + \frac1{ (n+m)^2 } \Big) + \\
+ \( \Im \hat b(2\pi m/T) - \Im \hat b(2\pi n/T) \)^2 \cdot \frac1{
  (n-m)^2 } + \\
+ \( \Im \hat b(2\pi m/T) + \Im \hat b(2\pi n/T) \)^2 \cdot \frac1{
  (n+m)^2 } \bigg) < \infty \, ;
\end{multline*}
it is enough to prove that
\begin{align}\label{9.a}
\sum_{m,n:\,0\le m<n} \frac{ \ln^{\al-1} (m+2) \ln^{\al-1} (n+2) }{
  (n-m)^2 } | \hat b (2\pi m/T) - \hat b (2\pi n/T) |^2 < \infty \, ,
  \\ \label{9.b}
\sum_{m,n:\,0\le m<n} \frac{ \ln^{\al-1} (m+2) \ln^{\al-1} (n+2) }{
  (n+m)^2 } \( \Im \hat b (2\pi m/T) + \Im \hat b (2\pi n/T) \)^2 <
  \infty \, .
\end{align}
We treat separately two cases, $ 0 \le m < \sqrt n $ and $ \sqrt n \le
m < n $. The first case, $ 0 \le m < \sqrt n $, is simple; just using
boundedness of $ \hat b $ we have for \eqref{9.a} and \eqref{9.b} as
well,
\begin{multline*}
\sum_{m,n:\,0\le m<\sqrt n} \dotsc \le \const \cdot \sum_{m,n:0\le
  m<\sqrt n} \frac{ \ln^{2\al-2} (n+2) }{ n^2 } \le \\
\le \const \cdot \sum_n \sqrt n \cdot \frac{ \ln^{2\al-2} (n+2) }{ n^2
  } < \infty \, .
\end{multline*}
We turn to the other case, $ \sqrt n \le m < n $. Now $ \ln(n+2) = O (
\ln(m+2) ) $. Lemma \ref{9.5} gives $ \Im \hat b(\la) = O \(
\frac1{\ln^\al \la} \) $, hence
\begin{multline*}
\frac{ \ln^{\al-1} (m+2) \ln^{\al-1} (n+2) }{ (n+m)^2 } \( \Im \hat b
  (2\pi m/T) + \Im \hat b (2\pi n/T) \)^2 = \\
= O \bigg( \frac{ \ln^{\al-1} (m+2) \ln^{\al-1} (n+2) }{ n^2 } \Big(
  \frac1{ \ln^{2\al} (m+2) } + \frac1{ \ln^{2\al} (n+2) } \Big) \bigg)
  = \\
= O \bigg( \frac{ \ln^{2\al-2} (n+2) }{ n^2 \ln^{2\al} (m+2) } \bigg)
  = O \bigg( \frac1{ n^2 \ln^2 (n+2) } \bigg) \, ;
\end{multline*}
summing over $ m $ gives $ O \( \frac1{ n \ln^2 (n+2) } \) $, a
convergent series in $ n $, which proves \eqref{9.b}.

It remains to prove the most delicate case, \eqref{9.a} for $ \sqrt n
\le m < n $. Neglecting a finite number of terms, we get $ m $ large
enough for using the asymptotic relation of Lemma \ref{9.5}:
\[
\hat b (2\pi m/T) - \hat b (2\pi n/T) = O \bigg( \int_{2\pi m/T}^{2\pi
n/T} \frac{ d\la }{ \la \ln^\al \la } \bigg) = O \bigg( \frac1{
\ln^{\al-1} m } - \frac1{ \ln^{\al-1} n } \bigg) \, .
\]
However, $ (\ln m)^{-(\al-1)} - (\ln n)^{-(\al-1)} \le (\al-1) (\ln
m)^{-\al} (\ln n - \ln m) $ and $ (\ln m)^{-\al} \le ( \frac12 \ln
n)^{-\al} $, therefore
\[
\hat b (2\pi m/T) - \hat b (2\pi n/T) = O \bigg( \frac{ \ln n - \ln m
}{ \ln^\al n } \bigg) \, .
\]
It is enough to prove that
\[
\sum_{m,n:\,\sqrt n \le m<n} \frac{ \ln^{\al-1} m \cdot \ln^{\al-1} n }{
(n-m)^2 } \bigg( \frac{ \ln n - \ln m }{ \ln^\al n } \bigg)^2 < \infty
\]
or, equivalently,
\[
\sum_n \frac1{ n^2 \ln^2 n } \sum_{m:\,\sqrt n\le m<n} \bigg( \frac{ \ln
\frac n m }{ 1 - \frac m n } \bigg)^2 < \infty \, .
\]
It remains to note that
\[
\frac1n \sum_{m:\,\sqrt n\le m<n} \bigg( \frac{ \ln \frac n m }{ 1 -
\frac m n } \bigg)^2 \le \int_0^1 \bigg( \frac{\ln(1/u)}{1-u} \bigg)^2
\, du < \infty \, .
\]
\end{proof}

\begin{proposition}\label{9.9}
If $ B $ satisfies \eqref{9.2} then $ X_k $ (defined by \eqref{9.3})
are orthogonal w.r.t.\ some admissible norm on the FHS-space $ G_{0,T}
$.
\end{proposition}

\begin{proof}
Here is an equivalent formulation: there is an operator $ A : l_2 \to
G_{0,T} $ such that
\[
A (c_1,c_2,\dotsc) = \sum_k \frac{ c_k }{ \|X_k\| } X_k \quad
\text{for all } (c_1,c_2,\dotsc) \in l_2 \, ,
\]
and $ A $ is an FHS-isomorphism, in other words, an equivalence
operator in the sense of Feldman \cite[Def.~1]{Fe}. It means that $ A
$ is one-to-one onto, has a bounded inverse, and $ \sqrt{A^*A} - I \in
\HS $ (the Hilbert-Schmidt class of operators). The latter is
equivalent to $ A^* A - I \in \HS $, see \cite[Lemma
1(b)]{Fe}.\footnote{%
 Though, his formulation of the lemma is incorrect, see the review
 21\#1546 in Mathematical Reviews.}

Matrix elements of $ A^* A $ are $ \dfrac{ | \ip{X_m}{X_n} | }{ \| X_m
\| \| X_n \| } $; Lemma \ref{9.6} shows that $ A^* A - I \in \HS $
and, of course, $ A $ is bounded. It remains to prove that $ A $ has a
bounded inverse. The range of $ A $ being evidently dense, we have to
prove that $ \| A x \| \ge \eps \| x \| $ for some $ \eps $, that is,
$ 0 $ does not belong to the spectrum of $ A^* A $. The spectrum
accumulates to $ 1 $ only (since $ A^* A - I \in \HS $); we have to
prove that $ 0 $ is not an eigenvalue, that is,
\[
\sum_k \frac{ c_k }{ \| X_k \| } X_k = 0 \imply c_1 = c_2 = \dots = 0
\]
for all $ (c_1,c_2,\dotsc) \in l_2 $. It is enough to prove that the
following formula is a correct definition of (continuous) linear
functionals $ X^1, X^2, \dotsc $ on $ G_{0,T} $:
\[
X^k (g) = \frac1T \int_0^T g(t) \exp (-2\pi i k t/T) \, dt \quad
\text{for } g \in G_{0,T} \, ;
\]
indeed, it will follow that
\[
c_k = \| X_k \| \cdot X^k \bigg( \sum_l \frac{ c_l }{ \| X_l \| } X_l
\bigg) \, .
\]
The norm on $ G_{0,T} $, defined in terms of $ B(\cdot) $, uses $ B(t)
$ for $ t \in [-T,T] $ only. Therefore we may assume that $ B(t) $
vanishes outside of some bounded interval (and still satisfies
\eqref{9.2}). For every $ g \in L_2 (0,T) \subset G_{0,T} $,
\[
\| g \|^2_{G_{0,T}} = \int_{-\infty}^{+\infty} \hat B(\la) | \hat
g(\la) |^2 \, d\la \, ;
\]
here $ \hat B $ is the Fourier transform (normalized as to be unitary)
of $ B $, and $ \hat g $ --- of $ g $. The formula $ (Zg) (\la) =
\sqrt{ \hat B(\la) } \hat g(\la) $ defines a linear isometric
embedding $ Z : G_{0,T} \to L_2 (\R) $ on the dense subset $ L_2 (0,T)
\subset G_{0,T} $; we may extend $ Z $ to the whole $ G_{0,T} $ by
continuity. Every $ \phi \in L_2 (\R) $ gives a linear functional on $
G_{0,T} $, namely, $ g \mapsto \int \sqrt{ \hat B(\la) } \hat g(\la)
\phi(\la) \, d\la $. In order to get $ X^k $, we take $ \phi $ such
that $ \sqrt{ \hat B(\la) } \phi (\la) $ is the Fourier transform of $
(1/T) \exp (-2\pi i k t/T) \One_{(0,T)} $; it remains to verify that
such $ \phi $ belongs to $ L_2 (\R) $.

The function $ (1/T) \exp (-2\pi i k t/T) \One_{(0,T)} $ is of finite
variation; its Fourier transform is $ O \( 1/\sqrt{1+\la^2} \) $; it
remains to check that
\[
\int \frac1{ (1+\la^2) \hat B (\la) } \, d\la < \infty \, .
\]
However, the continuous function $ \hat B $ never vanishes, and
$\displaystyle \hat B (\la) \sim \frac{\const}{ \ln^{\al-1} |\la| } $
for $ \la \to \pm\infty $ by Lemma \ref{9.5}.
\end{proof}

\begin{proposition}\label{9.10}
If $ B $ satisfies \eqref{9.2} then $ G_{-T,T} = G_{-T,0} \oplus
G_{0,T} $ (in the FHS sense).
\end{proposition}

\begin{proof}
Vectors $ X_k, Y_k $ (defined by \eqref{9.3}) are orthogonal w.r.t.\
some admissible norm on the FHS-space $ G_{-T,T} $; the proof is quite
similar to the proof of Proposition \ref{9.9}.
\end{proof}

Combining Proposition \ref{9.10} with elementary properties of spaces
$ G_{a,b} $ and operators $ U_t $ mentioned in the beginning of the
section, we get the following result.

\begin{theorem}\label{9.11}
If $ B $ satisfies \eqref{9.2} for some $ \al \in (1,\infty) $, then $
\( (G_{a,b}), (U_t) \) $ is a sum system (as defined by \ref{7.1}).
\end{theorem}

%% file: sect10.tex
Sum systems given by Theorem \ref{9.11} depend on the parameter $ \al
\in (1,\infty) $.\footnote{%
 That is, each such system corresponds to some $ \al $. However, for a
 given $ \al $ there is some freedom when choosing $ B $ satisfying
 \eqref{9.2}.}
Their exponentials, given by Theorem \ref{7.5}, are product
systems. Such product systems for different $ \al $ are nonisomorphic,
which will be shown using Proposition \ref{8.1}. Accordingly, we
consider a sequence of elementary sets $ E_n \subset (0,1) $, and we
want to know, which sequences $ (E_n) $ satisfy $ \liminf
G_{(0,1)\setminus E_n} = G_{(0,1)} $ and $ \limsup G_{E_n} = \{0\} $;
here $ G_{E_n} $ (as well as $ G_{(0,1)\setminus E_n} $) correspond
(as explained in Sect.~8) to the sum system given by Theorem
\ref{9.11}; recall that $ \liminf $ was defined by \ref{3.6}, and $
 \limsup $ by \ref{3.13}.

\begin{lemma}\label{10.1}
If $ \mes E_n \to 0 $ then $ \liminf G_{(0,1)\setminus E_n} =
G_{(0,1)} $. (Here $ \mes E_n $ is Lebesgue measure of $ E_n $.)
\end{lemma}

\begin{proof}
We have to represent an arbitrary vector $ g \in G_{(0,1)} $ as $ \lim
g_n $ for some $ g_n \in G_{(0,1)\setminus E_n} $. It is enough to
consider a dense set of vectors $ g $ (since $ \liminf $ is always
closed). Let $ g \in L_2 (0,1) \subset G_{(0,1)} $ and $ g_n = g \cdot
\One_{(0,1)\setminus E_n} \in G_{(0,1)\setminus E_n} \cap L_2(0,1)
$. Clearly, $ \mes E_n \to 0 $ implies $ g \cdot \One_{E_n} \to 0 $ in
$ L_2(0,1) $, hence $ g_n \to g $ in $ L_2(0,1) $, therefore $ g_n \to
g $ in $ G_{0,1} $.
\end{proof}

\begin{lemma}\label{10.2}
Assume that $ B $ satisfies \eqref{9.2} for a given $ \al \in
(1,\infty) $, and $ \mes E_n = o ( 1 / \ln^{\al-1} n ) $, and $ E_n $
consists of (no more than) $ n $ intervals. Then $ \limsup G_{E_n} =
\{ 0 \} $.
\end{lemma}

\begin{proof}
Let $ g_n \in G_{E_n} $, $ \| g_n \| \le 1 $; we have to prove that $
g_n \to 0 $ weakly. Introduce
\[
h_n (t) = c \frac{ n^2 }{ \ln^{(\al-1)/2} n } \cdot
\One_{(-1/n^2,+1/n^2)}
\]
where $ c > 0 $ is chosen such that for all $ n $ large enough,
\[
B - h_n * h_n \quad \text{is positively definite;}
\]
in terms of Fourier transform it means that
\[
\( \hat h_n (\la) \)^2 \le B(\la) \quad \text{for all } \la \in \R \,
;
\]
such $ c $ exists due to the asymptotic relation $ \hat B (\la) \sim
\frac{\const}{ \ln^{\al-1} |\la| } $ for large $ |\la| $ (recall
\ref{9.5}). The positive definiteness means that
\[
\| g * h_n \|_{L_2(\R)} \le \| g \|_{G_{0,1}}
\]
for all $ g \in L_2 (0,1) \subset G_{0,1} $ and then (extending the
convolution operator by continuity) for all $ g \in G_{0,1} $.

Denote by $ E'_n $ the $ (1/n^2) $-neighborhood of $ E_n $, then $
\mes E'_n \le (2n/n^2) + \mes E_n = o ( 1 / \ln^{\al-1} n ) $, and $ g
* h_n \in L_2 (E'_n) \subset L_2(\R) $.

We have to prove that $ \phi (g_n) \to 0 $ for every linear functional
$ \phi $ on $ G_{0,1} $. We may restrict ourselves to a dense subset
of functionals $ \phi $ (since $ \| g_n \| \le 1 $). In particular, we
may consider only functionals $ \phi_k $ defined by
\[
\phi_k (g) = \int g(t) \exp (2\pi i k t) \, dt \, .
\]
Taking into account that $\displaystyle \phi_k ( g * h_n ) = \frac{ 2c
}{ \ln^{(\al-1)/2} n } \cdot \frac{ \sin (2\pi k/n^2) }{ 2\pi k/n^2 }
\phi_k (g) $ we see that the following would be enough:
\[
\frac1{2c} \( \ln^{(\al-1)/2} n \) \phi_k ( g_n * h_n ) \to 0 \quad
\text{when } n \to \infty
\]
for every $ k $.

Recalling that $ g_n * h_n \in L_2 (E'_n) $ and $ \| g_n * h_n
\|_{L_2(\R)} \le \| g_n \|_{G_{0,1}} \le 1 $ we have
\[
| \phi_k ( g_n * h_n ) | \le \sqrt{ \mes E'_n } \cdot \| g_n * h_n
  \|_{L_2(\R)} = o \Big( \frac1{ \ln^{(\al-1)/2} n } \Big) \, .
\]
\end{proof}

\begin{lemma}\label{10.3}
Assume that $ B $ satisfies \eqref{9.2} for a given $ \al \in
(1,\infty) $. Then there exists a sequence of elementary sets $ E_n
\subset (0,1) $ such that $ \mes E_n = O ( 1 / \ln^{\al-1} n ) $, and
$ E_n $ consists of $ n $ intervals, and the relation $ \limsup
G_{E_n} = \{ 0 \} $ is violated.
\end{lemma}

\begin{proof}
The construction is straightforward, just $ n $ equidistant intervals
of equal length:
\[\begin{split}
E_n &= \bigcup_{k=1}^n \bigg( \frac1n \Big( k - \frac12 - \frac1{ n
  \ln^{\al-1} n } \Big), \frac1n \Big( k - \frac12 + \frac1{ n
  \ln^{\al-1} n } \Big) \bigg) \, ; \\
\mes E_n &= \frac2{ \ln^{\al-1} n } \, .
\end{split}\]
We have to find $ g_n \in G_{E_n} $ and $ g \in G_{0,1} $ such that
\begin{gather}
\label{10.a}
\sup_n \| g_n \|_{G_{0,1}} < \infty \, , \\
\label{10.b}
\limsup_n | \ip{ g_n }{ g }_{G_{0,1}} | > 0 \, .
\end{gather}
Still, the construction is straightforward:
\[\begin{split}
g_n &= ( 1 / \mes E_n ) \cdot \One_{E_n} = \frac12 \ln^{\al-1} n \cdot
  \One_{E_n} \, , \\
g &= \One_{(0,1)} \, ;
\end{split}\]
the proof of \eqref{10.a}, \eqref{10.b} is more complicated. Fourier
transform $ f \mapsto \hat f $ will be used, $ \hat f (\la) = \int
e^{i\la t} f(t) \, dt $. Recall that
\[
2\pi \ip{ f_1 }{ f_2 }_{G_{0,1}} = \int_{-\infty}^{+\infty} \hat B
(\la) \hat f_1 (\la) \overline{ \hat f_2 (\la) } \, d\la
\]
for all $ f_1, f_2 \in G_{0,1} $, and $ \hat B(\la) \sim
\frac{ \const }{ \ln^{\al-1} |\la| } $ for large $ |\la| $. An
elementary calculation gives
\[\begin{split}
\hat g (\la) &= \frac{ e^{i\la} - 1 }{ i\la } \, ; \\
& \sum_{k=0}^{n-1} \exp (ik\la/n) = \frac{ 1 - \exp(i\la) }{ 1 -
  \exp(i\la/n) } \, ; \\
\hat g_n (\la) &= (\ln n)^{\al-1} \exp \Big( i\frac{\la}{2n} \Big)
  \frac{ 1 - \exp(i\la) }{ 1 - \exp(i\la/n) } \frac1\la \sin \frac{
  \la }{ n \ln^{\al-1} n } \, ; \\
| \hat g_n (\la) | &= (\ln n)^{\al-1} \frac{ | \sin \frac\la2 | }{ |
  \sin \frac{\la}{2n} | } \frac1{|\la|} \Big| \sin \frac{ \la }{ n
  \ln^{\al-1} n } \Big| \, .
\end{split}\]
Note that
\[
\frac1{2\pi n} \int_0^{2\pi n} \frac{ \sin^2 \frac\la2 }{ \sin^2
\frac{\la}{2n} } \, d\la = \frac1{2\pi n} \int_0^{2\pi n} \Big|
\sum_{k=0}^{n-1} \exp (ik\la/n) \Big|^2 \, d\la = n \, .
\]
We have
\[
2\pi \| g_n \|^2 = \int_{-\infty}^{+\infty} \hat B(\la) | \hat g_n
(\la) |^2 \, d\la = 2 \bigg( \int_0^n + \int_n^\infty \bigg) \hat
B(\la) | \hat g_n (\la) |^2 \, d\la \, .
\]

For $ \la \in (0,n) $ we note that
\[\begin{split}
& \hat B (\la) \le \hat B(0) \, , \\
& | \hat g_n (\la) | \le \frac1n \frac{ | \sin \frac\la2 | }{ | \sin
  \frac{\la}{2n} | } \, , \\
& \int_0^n | \hat g_n (\la) |^2 \, d\la \le \frac1{n^2} \int_0^n
\frac{ \sin^2 \frac\la2 }{ \sin^2 \frac\la{2n} } \, d\la \le 2\pi \, ,
\end{split}\]
therefore
\[
\int_0^n \hat B (\la) | \hat g_n (\la) |^2 \, d\la \le 2\pi \hat B (0)
\quad \text{for all } n \, .
\]
For $ \la \in (n,\infty) $ we note that $ \hat B (\la) \le \frac{
\const }{ \ln^{\al-1} n } $ and change the scale, introducing
$\displaystyle u = \frac{ \la }{ n \ln^{\al-1} n } $:
\begin{multline*}
\int_n^\infty \hat B (\la) | \hat g_n (\la) |^2 \, d\la \le \frac{
  \const }{ \ln^{\al-1} n } ( \ln n )^{2(\al-1)} \int_n^\infty \frac{
  \sin^2 \frac\la2 }{ \sin^2 \frac{\la}{2n} } \frac1{\la^2} \sin^2
  \frac{ \la }{ n \ln^{\al-1} n } \, d\la \le \\
\le \frac{ \const }{ n } \int_0^\infty \frac{ \sin^2 (n \om_n u ) }{
  \sin^2 ( \om_n u ) } \frac{ \sin^2 u }{ u^2 } \, du \, ,
\end{multline*}
where $ \om_n = \frac12 \ln^{\al-1} n \to \infty $. On each period, $
u \in \( \frac{\pi}{\om_n} k, \frac{\pi}{\om_n}(k+1) \) $, we
substitute $ \frac{ \sin^2 u }{ u^2 } $ by its maximal value:
\begin{multline*}
\int_0^\infty \frac{ \sin^2 ( n \om_n u ) }{ \sin^2 ( \om_n u ) }
  \frac{ \sin^2 u }{ u^2 } \, du \le \\
\le ( 1 + o(1) ) \bigg( \int_0^\infty
  \frac{ \sin^2 u }{ u^2 } \, du \bigg) \underbrace{ \bigg(
  \frac{\om_n}{\pi} \int_0^{\pi/\om_n} \frac{ \sin^2 ( n \om_n u ) }{
  \sin^2 ( \om_n u ) } \, du \bigg) }_{=n} \, .
\end{multline*}
Hence $ \int_n^\infty \hat B (\la) | \hat g_n (\la) |^2 \, d\la = O(1)
$. So, \eqref{10.a} is verified.

Further,
\[
2\pi \ip{ g_n }{ g } = \int_{-\infty}^{+\infty} \hat B (\la) \hat g_n
(\la) \overline{ \hat g (\la) } \, d\la = \bigg( \int_{|\la|<M} +
\int_{|\la|>M} \bigg) \hat B (\la) \hat g_n (\la) \overline{ \hat g
(\la) } \, d\la
\]
for any $ M \in (0,\infty) $. It is easy to see that
\[
\hat g_n (\la) \xrightarrow[n\to\infty]{} \hat g (\la) \quad
\text{uniformly in } \la \in [-M,M] \, ,
\]
which implies
\[
\int_{|\la|<M} \hat B (\la) \hat g_n (\la) \overline{ \hat g (\la) }
\, d\la \to \int_{|\la|<M} \hat B (\la) | \hat g (\la) |^2 \, d\la
\quad \text{for } n \to \infty \, .
\]
On the other hand,
\begin{multline*}
\bigg| \int_{|\la|>M} \hat B (\la) \hat g_n (\la) \overline{ \hat g
  (\la) } \, d\la \bigg| \le \\
\le \underbrace{ \bigg( \int_{|\la|>M} \hat B
  (\la) | \hat g_n (\la) |^2 \, d\la \bigg)^{1/2} }_{\le\sup_n \|g_n\|}
  \underbrace{ \bigg( \int_{|\la|>M} \hat B (\la) | \hat g (\la) |^2 \,
  d\la \bigg)^{1/2} }_{\to 0 \text{ for } M\to\infty} \, ,
\end{multline*}
it tends to $ 0 $ uniformly in $ n $, when $ M \to \infty $. Choose $
M $ and $ \eps > 0 $ such that
\[
\sup_n \bigg| \int_{|\la|>M} \hat B (\la) \hat g_n (\la) \overline{
\hat g (\la) } \, d\la \bigg| \le \frac12 \eps \quad \text{and} \quad
\int_{|\la|<M} \hat B (\la) | \hat g (\la) |^2 \, d\la \ge \eps \, ,
\]
then
\[
\liminf_n \Re \ip{ g_n }{ g } \ge \frac12 \frac1{2\pi} \eps \, ,
\]
which implies \eqref{10.b}.
\end{proof}

\begin{theorem}
Let\footnote{%
 Here $ B' $ does not mean the derivative of $ B $, sorry.}
$ B', B'' $ satisfy \eqref{9.2} for some $ \al', \al'' $ respectively,
$ \al', \al'' \in (1,\infty) $, $ \al' \ne \al'' $. Then the
 corresponding product systems are nonisomorphic.
\end{theorem}

\begin{proof}
Suppose that $ \al' < \al'' $. Lemma \ref{10.3} gives elementary sets
$ E_n \subset (0,1) $ such that $ \mes E_n = O ( 1 / \ln^{\al''-1} n )
$, and $ E_n $ consists of $ n $ intervals, and the relation $ \limsup
G''_{E_n} = \{ 0 \} $ is violated. Taking into account that $ O ( 1 /
\ln^{\al''-1} n ) = o ( 1 / \ln^{\al'-1} n ) $ we get $ \limsup
G'_{E_n} = \{ 0 \} $ by Lemma \ref{10.2}. Also, $ \liminf
G'_{(0,1)\setminus E_n} = G'_{(0,1)} $ and $ \liminf
G''_{(0,1)\setminus E_n} = G''_{(0,1)} $ by Lemma \ref{10.1}. So, $
E_n $ satisfy \ref{8.1}(a) but violate \ref{8.1}(b). By Proposition
\ref{8.1} the product systems are nonisomorphic.
\end{proof}

\begin{lemma}\label{10.7}
Assume that $ B $ satisfies \eqref{9.2} for a given $ \al \in
(1,\infty) $, and
\[\begin{split}
& E_n = \bigcup_{k=1}^n \bigg( \frac1n \Big( k - \frac12 - \frac{ \mes
  E_n }{ 2n } \Big), \frac1n \Big( k - \frac12 + \frac{ \mes E_n }{ 2n
  } \Big) \bigg) \, , \\
& (\ln n)^{\al-1} \mes E_n \to \infty \quad \text{for } n \to \infty \,
  , \\
& g_n = ( 1 / \mes E_n ) \cdot \One_{E_n} \, , \\
& g = \One_{(0,1)} \, .
\end{split}\]
Then $ \| g_n - g \|_{G_{0,1}} \to 0 $ for $ n \to \infty $.
\end{lemma}

\begin{proof}
Similarly to the proof of Lemma \ref{10.3} we have for any $ M \in
(0,\infty) $
\[
\hat g_n (\la) \xrightarrow[n\to\infty]{} \hat g (\la) \quad
\text{uniformly on } \la \in [-M,M] \, .
\]
This time, however, the following compactness property holds:
\[
\int_{|\la|>M} \hat B (\la) | \hat g_n (\la) |^2 \, d\la
\xrightarrow[M\to\infty]{} 0 \quad \text{uniformly in } n \, ,
\]
which ensures $ \| g_n - g \| \to 0 $. In order to prove the
compactness property we estimate integrals similarly to the proof of
\ref{10.3}. We have
\[
| \hat g_n (\la) | = \frac2{\mes E_n} \frac{ | \sin \frac\la2 | }{ |
\sin \frac\la{2n} | } \frac1{|\la|} \Big| \sin \frac{ \la \mes E_n }{
2n } \Big| \, .
\]
For $ \la \in (M,n) $ we note that\footnote{%
 Assuming that $ M $ is large enough.}
\[\begin{split}
& \hat B (\la) \le \frac{ \const }{ \ln^{\al-1} M } \, , \\
& | \hat g_n (\la) | \le \frac1n \frac{ | \sin \frac\la2 | }{ |
  \sin \frac\la{2n} | } \, , \qquad \text{(the same as in \ref{10.3})}
  \\
& \int_M^n \hat B (\la) | \hat g_n (\la) |^2 \, d\la \le \frac{ \const
  }{ \ln^{\al-1} M } \int_0^n | \hat g_n (\la) |^2 \, d\la \le \frac{
  \const }{ \ln^{\al-1} M } \cdot 2\pi \xrightarrow[M\to\infty]{} 0 \,
  .
\end{split}\]
For $ \la \in (n,\infty) $, introducing $ \displaystyle u = \frac{
\mes E_n }{ 2n } \la $ and $ \om_n = 1 / \mes E_n \to \infty $, we
have
\begin{multline*}
\int_n^\infty \hat B (\la) | \hat g_n (\la) |^2 \, d\la \le \frac{
  \const }{ \ln^{\al-1} n } \bigg( \frac2{ \mes E_n } \bigg)^2
  \int_n^\infty \frac{ \sin^2 \frac\la2 }{ \sin^2 \frac\la{2n} }
  \frac1{\la^2} \sin^2 \frac{ \la \mes E_n }{ 2n } \, d\la \le \\
\le \underbrace{ \frac{ \const }{ \ln^{\al-1} n } \frac2{ \mes E_n }
  }_{\to0} \cdot \frac1n \cdot \underbrace{ \int_0^\infty \frac{
  \sin^2 ( n \om_n u ) }{ \sin^2 ( \om_n u ) } \frac{ \sin^2 u }{ u^2
  } \, du }_{O(n)} \xrightarrow[n\to\infty]{} 0 \, .
\end{multline*}
\end{proof}

\begin{proposition}
If $ B $ satisfies \eqref{9.2} for some $ \al \in (1,\infty) $ then
the corresponding product system is unitless.
\end{proposition}

\begin{proof}
Assume the contrary, then there exist $ \psi_{a,b} \in H_{a,b} = \Exp
G_{a,b} $ such that $ \| \psi_{a,b} \| = 1 $ and $ \psi_{a,b} \otimes
\psi_{b,c} = \psi_{a,c} $ whenever $ a < b < c $. We'll use $
\psi_{a,b} $ for $ 0 \le a < b \le 1 $ only. Introduce a Gaussian type
space $ (\Om,\F,\P,G) $ such that $ G_{0,1} = G / \Const $; we have $
H_{a,b} = L_2 (\Om,\F_{a,b},\P|_{\F_{a,b}}) $ where $ \F_{a,b} \subset
\F $ is the \sif\ generated by the Gaussian subspace $ G_{a,b} \subset
G_{0,1} $. Consider measures $ \mu_{a,b} = |\psi_{a,b}|^2 $; as was
explained in Sect.~1, $ \mu_{a,b} $ is a probability measure on the
\sif\ $ \F_{a,b} $ such that
\[
\frac{ \mu_{a,b} }{ P|_{\F_{a,b}} } = \bigg| \frac{ \psi_{a,b} }{
\sqrt{ P|_{\F_{a,b}} } } \bigg|^2 \quad \text{for all } P \in \P \, .
\]
The relation $ \psi_{a,b} \otimes \psi_{b,c} = \psi_{a,c} $ implies
\[
\mu_{a,b} \otimes \mu_{b,c} = \mu_{a,c} \, .
\]
In other words, \sif s $ \F_{a,b} $ and $ \F_{b,c} $ are independent
w.r.t.\ the measure $ \mu = \mu_{0,1} $ on $ \F $, and $ \mu_{a,b} $
is just the restriction of $ \mu $ to $ \F_{a,b} $.\footnote{%
 Clearly, $ \mu $ is a probability measure on $ (\Om,\F) $, absolutely
 continuous w.r.t.\ any $ P \in \P $. However, $ P $ need not be
 absolutely continuous w.r.t.\ $ \mu $; that is, $ \mu $ need not
 belong to $ \P $.}
Moreover, \sif s $ \F_{t_0,t_1}, \F_{t_1,t_2}, \dots, \F_{t_{n-1},t_n}
$ are $ \mu $-independent whenever $ 0 \le t_0 < t_1 < \dots < t_n \le
1 $. Every elementary set $ E \subset (0,1) $ determines its sub-\sif\
$ \F_E \subset \F $; as explained in Sect.~8,
\[
E = (t_0,t_1) \cup (t_2,t_3) \cup \dots \cup (t_{2n},t_{2n+1}) \imply
\F_E = \F_{t_0,t_1} \otimes \dots \otimes \F_{t_{2n},t_{2n+1}}
\]
whenever $ 0 \le t_0 < t_1 < \dots < t_{2n+1} \le 1 $. Note that $
\F_{E_1}, \F_{E_2} $ are $ \mu $-independent whenever $ E_1 \cap E_2 =
\emptyset $.

Consider elementary sets
\[\begin{split}
E_n &= \Big( 0, \frac1{2n} \Big) \cup \Big( \frac2{2n},\frac3{2n}
  \Big) \cup \dots \cup \Big( \frac{2n-2}{2n}, \frac{2n-1}{2n} \Big)
  \, , \\
\mes E_n &= \frac12 \, ,
\end{split}\]
and vectors
\[\begin{split}
g_n &= 2 \cdot \One_{E_n} \in G_{E_n} \subset G_{0,1} \, , \\
g &= \One_{(0,1)} \in G_{0,1} \, , \\
h_n &= 2 \cdot \One_{(0,1)\setminus E_n} \in G_{(0,1)\setminus E_n}
  \subset G_{0,1} \, .
\end{split}\]
Lemma \ref{10.7} shows\footnote{%
 Shifting the set $ E_n $ of Lemma \ref{10.7} does not invalidate the
 lemma.}
that $ g_n \to g $ and $ h_n \to g $ in $ G_{0,1} $. Though, $ G_{0,1}
$ is not a subspace but a quotient space $ G / \Const $ of $ G $;
anyway, we may choose elements of $ G $, denoted again by $ g_n, g,
h_n $, such that
\[\begin{split}
& g_n \in L_0 (\Om,\F_{E_n},\P) \, , \\
& g \in L_0 (\Om,\F,\P) \, , \\
& h_n \in L_0 (\Om,\F_{(0,1)\setminus E_n},\P) \, , \\
& g_n \to g \quad \text{and} \quad h_n \to g \quad \text{in } L_0
(\Om,\F,\P) \, .
\end{split}\]
The natural map\footnote{%
 Generally, non-invertible, since $ \mu $ need not belong to $ \P $.}
$ L_0 (\Om,\F,\P) \to L_0 (\Om,\F,\mu) $ allows us to treat $ g_n, g,
h_n $ as elements of $ L_0 (\Om,\F,\mu) $. Now they are random
variables; $ g_n \to g $ and $ h_n \to g $ in probability. On the
other hand, for every $ n $, the two random variables $ g_n, h_n $ are
independent (since $ \F_{E_n} $ and $ \F_{(0,1)\setminus E_n} $ are $
\mu $-independent). It follows that $ g $ is independent of itself,
that is, $ g $ is constant $ \mu $-almost sure.\footnote{%
 Proof: let $ \phi : \R \to \R $ be a continuous and bounded function,
 then $ \int \phi( g_n(\om) ) \phi( h_n(\om) ) \, \mu(d\om) = \( \int
 \phi( g_n(\om) ) \, \mu(d\om) \) \cdot \( \int \phi( h_n(\om) ) \,
 \mu(d\om) \) $ due to independence. The limit for $ n \to \infty $
 gives $ \int \( \phi( g(\om) ) \)^2 \, \mu(d\om) = \( \int \phi(
 g(\om) ) \, \mu(d\om) \)^2 $, which is impossible unless $ g = \const
 $.}

Consider a Gaussian measure $ \ga \in \P_G $. Though, $ g $ need not
be constant $ \ga $-almost sure; however, $ g $ must be constant on a
set of positive probability w.r.t.\ $ \ga $. On the other hand, the
distribution of $ g $ w.r.t.\ $ \ga $ is normal (Gaussian); it cannot
have an atom unless it is degenerate, which means that $ \| g
\|_{G_{0,1}} $ must vanish. However, it does not vanish, which is
evident when using Fourier transform. A contradiction.
\end{proof}